\documentclass[11pt]{article}

\usepackage[sort]{cite}       % nicer citations
\usepackage{amssymb}          % various ams stuff
\usepackage{latexsym}         % dito
\usepackage{amsmath}          % dito
\usepackage{amscd}            % nasty commutative ams diagrams
\usepackage[latin1]{inputenc} % smart input of special characters
\usepackage{eucal}            % nicer caligraphic fonts
\usepackage{bbm}              % nicer bbm fonts
\usepackage{theorem}          % nettere Theoreme
\usepackage{graphicx}         % neuere Art, eps-Bildchen einzubinden
\usepackage{stmaryrd}         % noch ein paar mehr Symbole (Klammern)
\usepackage{mathrsfs}         % nette kaligraphische Schrift
\usepackage[ps,dvips]{xypic}  % the cool diagrams need this with ps
\usepackage{diagxy}           % the real diagrams !

\title{Locally Noncommutative Space-Times}

\author{\textbf{Dorothea Bahns}\thanks{Department f\"ur Mathematik,
  Universit\"at Hamburg,  Bundesstr. 55,
  D-20253 Hamburg,  Germany, 
{\tt   bahns@math.uni-hamburg.de}} \addtocounter{footnote}{2}
  and
  \textbf{Stefan Waldmann}\thanks{  Fakult\"at f\"ur Mathematik und Physik,
  Albert-Ludwigs-Universit{\"a}t Freiburg,
  Physikalisches Institut,  Hermann-Herder-Stra{\ss}e 3,
  D 79104 Freiburg,
  Germany,
{\tt    Stefan.Waldmann@physik.uni-freiburg.de}
}
}

\date{July 2006\\[0.5cm] FR-THEP 2006/12, ZMP-HH/2006-13}

\renewcommand{\mathbb}[1]{\mathbbm{#1}} % use nicer bbm fonts

\newcounter{comment}

\newcommand{\Lie}        {\operatorname{\mathscr{L}\!}}    
\newcommand{\cc}[1]      {\overline{{#1}}}              
\newcommand{\id}         {\operatorname{\mathsf{id}}}   
\newcommand{\supp}       {\operatorname{\mathrm{supp}}}  
   
\newcommand{\tr}         {\operatorname{\mathsf{tr}}}

\newcommand{\End}        {\operatorname{\mathsf{End}}}   
\newcommand{\SP}[1]      {\left\langle{#1}\right\rangle} 
\newcommand{\Unit}       {\mathbb{1}}

\newcommand{\I}          {\mathrm{i}}
\newcommand{\E}          {\mathrm{e}}
\newcommand{\D}          {\operatorname{\mathrm{d}}}

\newcommand{\pr}         {\mathrm{pr}}

\newcommand{\Var}      	{\operatorname{\mathrm{Var}}}

\newcommand{\tsf}	{\textstyle \frac}
\newcommand{\R} 	{\mathbb R}

\newcommand{\Ver}        {\operatorname{\mathrm{Ver}}}

\newcommand{\Anti}       {\Lambda}
\newcommand{\Sym}        {\mathrm{S}}
\newcommand{\Xver}       {\operatorname{\mathfrak{X}}_{\scriptscriptstyle{\mathrm{ver}}}}
\newcommand{\ver}        {\mathsf{v}}
\newcommand{\Schouten}[1]{\left\llbracket{#1}\right\rrbracket}

\newcommand{\HCdiff}   {\mathrm{HC}_{\scriptscriptstyle{\mathrm{diff}}}}
\newcommand{\HHdiff}   {\mathrm{HH}_{\scriptscriptstyle{\mathrm{diff}}}}
\newcommand{\HCdiffver}{\mathrm{HC}_{\scriptscriptstyle{\mathrm{diff,ver}}}}
\newcommand{\HHdiffver}{\mathrm{HH}_{\scriptscriptstyle{\mathrm{diff,ver}}}}
\newcommand{\Uver}     {\operatorname{\mathrm{U}}_{\scriptscriptstyle{\mathrm{ver}}}}
\newcommand{\uUver}    {\operatorname{\underline{\mathrm{U}}}_{\scriptscriptstyle{\mathrm{ver}}}}

\newcommand{\starp}      {\mathbin{\star_p}}
\newcommand{\tstar}      {\mathbin{\tilde{\star}}}
\newcommand{\tstarp}     {\mathbin{\tilde{\star}_p}}

\newtheorem{lemma} {Lemma} [section]
\newtheorem{proposition} [lemma] {Proposition}

\newtheorem{theorem} [lemma] {Theorem}

\theorembodyfont{\rmfamily}
\newtheorem{example}[lemma]{Example}
\newtheorem{remark}[lemma]{Remark}

\newenvironment{proof}[1][{}]{
  \par\noindent
  \textsc{Proof{#1}:}
}
{
  \hspace*{\fill} $\blacksquare$\newline
}

\numberwithin{equation}{section}

\newcommand{\beqa}	{\begin{eqnarray*}}
\newcommand{\eeqa}	{\end{eqnarray*}}

\newcommand{\beqan}	{\begin{eqnarray}}
\newcommand{\eeqan}	{\end{eqnarray}}

\begin{document}

\maketitle

\begin{abstract}
    \noindent
    Localized noncommutative structures for manifolds with connection
    are constructed based on the use of vertical star products. The
    model's main feature is that two points that are far away from
    each other will not be subject to a deviation from classical
    geometry while space-time becomes noncommutative for pairs of
    points that are close to one another.
\end{abstract}

\section{Introduction}
\label{sec:Intro}

Models of noncommutative space-times have become increasingly popular
in the recent past and are believed to be reasonable candidates for
Planck scale physics, at least as an intermediate step towards a full
quantum mechanical treatment of geometry. The main idea is that in the
regime where quantum theory and general relativity are no longer
independent, the notion of a point in space-time becomes meaningless,
and a finite minimal length or uncertainty relations for coordinate
functions have to be postulated in order to prohibit the localization
of points with arbitrarily high precision.  Models with uncertainty
relations are usually implemented by considering noncommuting
coordinate operators, replacing the coordinate functions.

More concretely, the situation can be seen as follows: for the
noncommutative algebra describing the noncommutative space-time one
chooses a symbol calculus resulting in a Poisson structure $\theta$ on
the classical space-time together with a corresponding quantization of
some reasonable class of functions on $M$. Quantization is achieved
either analytically, see e.g.  \cite{rieffel:1997a,
  doplicher.fredenhagen.roberts:1995a}, or in the context of
deformation quantization \cite{bayen.et.al:1978a} in form of a star
product, see e.g.~\cite{jurco.schupp.wess:2000a}.  A particularly
simple and well-studied example is the noncommutative Minkowski space
$M = \mathbb{R}^{1, 3}$ where the Poisson structure is chosen to be
\emph{constant} and symplectic. Of course, this is a highly
non-geometric situation and in view of an aim towards general
relativity, this model has to be understood as a toy-model that should
be surpassed by a more geometric approach.

Having established a noncommutative space-time, it is of particular
interest to study dynamics in it, usually in form of a quantum field
theory defined on such a space-time; and much work has been done in
that direction.  However, various technical but more importantly, also
conceptual problems occur in this case, one of which we want to
address here: the global nature of the Poisson structure, in
particular in the case of \emph{constant} $\theta$ on Minkowski
space-time, necessarily results in effects that are visible at large
distances.

In quantum mechanics, this situation occurs in the sense that large
momenta certainly play a crucial role in quantum effects, i.e. quantum
mechanics is very well observable at large distances in \emph{phase
  space}. However, there is a polarization separating momenta and
coordinates from one another and their respective `magnitudes' are
comparable only after a suitable scaling involving characteristic
parameters of the system as well as Planck's constant~$\hbar$.
Clearly, with regard to the commutation relations, the case of
constant nondegenerate $\theta$ in noncommutative Minkowski space is
mathematically equivalent to quantum mechanics in a space of half the
dimension. In this case, however, all noncommuting operators are of
the same kind, and short distance/long distance effects are expected
to mix.

In any case, more sophisticated investigations seem to indicate that a
constant $\theta$ leads to macroscopic effects that (if the theory's
predictions are taken seriously) would have to be visible. A prominent
example is the mixing of ultraviolet and infrared divergences that
appears in {\em Euclidean} noncommutative field theories. A related
effect is the modification of the dispersion relation in field
theories with Lorentz signature, where the largest deviation from the
ordinary relation appears at {\em small} momenta, see
e.g.~\cite{liao.sibold:2002a}.

Moreover, a violation of microcausality at all scales of distances
results from the application of the highly nonlocal twisted
convolution product (i.e. the Weyl-Moyal product in its integral
form), see however the discussion of the cluster decomposition
property in \cite{bahns:2006a}. It should be noted that for Lie-type
Poisson structures, i.e.  linear Poisson structures on Minkowski
space-time, or quadratic Poisson structures as arising from a quantum
group approach, this situation is even worse.

This leads to the following natural question: can we modify the noncommutativity
in such a way that it decays in a reasonable way for large distances?
Considering a non-constant $\theta$ that vanishes in  (spatial) infinity does of
course not solve the problem: in such a scenario, the noncommutativity would
vary within the universe, such that some regions  would have `more' or `less'
noncommutativity than others, depending on the absolute position.

Instead, our suggestion is to  take the concept of \emph{distance} as a starting
point and  consider the distance of \emph{two} points. Hence, in our approach,
$M \times M$ is endowed with a noncommutative structure, instead of $M$ itself
alone. It is now straightforward to define noncommutativity only at small
distances: we simply consider a Poisson structure and a
corresponding star product on $M \times M$ that is nontrivial only close to the
diagonal $\Delta_M$ in $M \times M$ and zero or quickly decreasing 
away from $\Delta_M$.

The aim of this paper now is to set up a kinematical framework for
such types of noncommutativity and explore some of their features. The
whole approach has to be understood as a first step as we have not yet
investigated any form of dynamics on our \emph{locally noncommutative
  space-time}. This will be subject to future projects.

Our paper is organized as follows: In Section~\ref{sec:SmallDistances}
we first discuss the general framework of star products on $M \times
M$ localized close to the diagonal and use the exponential map of the
Levi-Civit\`a connection of the space-time manifold $M$ to pull-back
everything to the tangent bundle $TM$ of $M$. As crucial condition the
\emph{verticality} of the Poisson structure and the star product is
discussed in detail. In Section~\ref{sec:PropertiesVertical} we
investigate further properties of the vertical star product and show
how they can be used to endow every point of $M$ with a small
noncommutative neighborhood. Then we discuss in detail states and the
corresponding expectation values of our observables including in
particular the measurement of the (Lorentz) metric on $M$ itself.
Section~\ref{sec:NoncommutativeMinkowski} is devoted to the particular
example of flat Minkowski space-time. Though our approach is fully
geometric in general, this provides an important and simple example
which we investigate in detail. The last
Section~\ref{sec:QuestionOutlook} contains a discussion of further
open questions and possible extensions and limitations of the model. In
Appendix~\ref{sec:VerticalFormality} we have included a detailed technical
discussion of vertical Poisson structures, vertical star products and
the vertical formality theorem governing their existence and
classification.

%
% many thanks to...
%

\noindent
\textbf{Acknowledgement:} We would like to thank Klaus Fredenhagen,
Jakob Heller, Stefan Weiß, Julius Wess and Jochen Zahn for valuable
discussions and remarks.

%%%%%%%%%%%%%%%%%%%%%%%%%%%%%%%%%%%%%%%%%%%%%%%%%%%%%%%%%%%%%%%%%
% Noncommutative Structure at Small Distances
%

%
% Noncommutative Structure at Small Distances
%

\section{Noncommutative Structure at Small Distances}
\label{sec:SmallDistances}

%
% General
%

Starting point of our construction is a smooth $n$-dimensional manifold $M$
which allows for a Lorentz metric $g$. Note however, that our construction only
depends on the respective Levi-Civit\`a connection, such that we may (and
frequently will) consider Riemannian metrics as well. Our aim is to define a
noncommutative structure on $M$ whose effects are visible only at small
distances in the sense that a pair of points  $q,q^\prime\in M$ will `feel'
noncommutative effects only when they are in the close vicinity of the diagonal
$\Delta_M$ in $M \times M$. 

%
% Classical space-time prerequisites
%

\subsection{Classical space-time prerequisites}
\label{subsec:ClassicalSpaceTime}

Let us first recall some well-known constructions from differential geometry
we shall need later on. The Levi-Civit\`a connection $\nabla$ of $g$
determines the geodesic structure of $M$ and thereby the exponential map
$\exp$, which is defined on some open neighborhood of the zero section of the
tangent bundle $\pi: TM \longrightarrow M$. We choose once and for all such an
open neighborhood $\mathcal{U} \subseteq TM$ of the zero section with the
property that the map
\begin{equation}
    \label{eq:DiffeoPhi}
    \Phi:
    TM \supseteq \mathcal{U} \ni v_p
    \; \mapsto \;
    \Phi(v_p) = 
    \left(
        \exp_p\left(-v_p\right), 
        \exp_p\left(v_p\right)
    \right) 
    \in M \times M
    \quad 
    \textrm{where}
    \quad
    v_p \in T_pM,
\end{equation}
is a diffeomorphism onto its image $\mathcal{V} \subseteq M \times M$.
Here, $v_p \in T_pM$ denotes a tangent vector at $p \in M$ and
$\exp_p$ is the exponential map at $p$, i.e. $t \mapsto \exp_p(tv_p)$
is the geodesic through $p$ at $t=0$ with initial velocity $v_p$.
Clearly, $\Phi$ maps the zero section diffeomorphically to the
\emph{diagonal} $\Delta_M \subseteq M \times M$ and $\mathcal{V}$ is
an open neighborhood of $\Delta_M$.

For $p \in M$ we set $\mathcal{U}_p = \mathcal{U} \cap T_pM$ which is
an open neighborhood of $0_p \in T_pM$. Its image $\mathcal{V}_p =
\exp_p\left(\mathcal{U}_p\right)$ is an open neighborhood of $p \in
M$. As usual, we obtain the well-known \emph{normal coordinates} $x^1,
\ldots, x^n: \mathcal{V}_p \longrightarrow \mathbb{R}$ on
$\mathcal{V}_p$ by choosing \emph{linear} coordinates $\xi^1, \ldots,
\xi^n$ with respect to some vector space basis $e_1, \ldots, e_n$ on
$T_pM$ and setting $x^i = \xi^i \circ \exp_p^{-1}$.  As $T_pM$ is
equipped with the metric $g_p$ we can even choose $e_1, \ldots, e_n$
to be orthonormal.

The geometric interpretation of $\Phi$ is now the following: suppose
$(q, q') \in \mathcal{V} \subseteq M \times M$ are within the image of
$\Phi$. Then $\Phi^{-1}(q, q') = v_p\in \mathcal{U}_p \subseteq T_pM$
and $v_p$ is such that $\exp_p(- v_p) = q$ and $\exp_p(v_p) = q'$.
Thus the point $p$ is the \emph{geodesic midpoint} between $q$ and
$q'$. The normal coordinates around $p$ can be seen as the geodesic
relative coordinates of $(q, q')$ with respect to their 'center of
mass' $p$. This explains our definition of $\Phi$ since this way the
situation for $q$ and $q'$ becomes most symmetric: we denote by
$\tau_{M \times M}: M \times M \longrightarrow M \times M$ the global
diffeomorphism $\tau_{M \times M}(q, q') = (q', q)$ and by $\tau_{TM}:
TM \longrightarrow TM$ the global diffeomorphism $\tau_{TM}(v_p) =
-v_p$. Then we have
\begin{equation}
    \label{eq:PhitauInvariant}
    \Phi \circ \tau_{TM} = \tau_{M \times M} \circ \Phi
\end{equation}
on $\mathcal{U} \cap \tau_{TM}(\mathcal{U})$. Thus it is natural to
demand that $\mathcal{U}$ is invariant under the reflection
$\tau_{TM}$. We will always assume that this is the case. As a
consequence also $\mathcal{V}$ is symmetric under the exchange
$\tau_{M \times M}$.

%
% Noncommutativity at small distances
%

\subsection{Noncommutativity at small distances}
\label{subsec:NCSmall}

We will now give a definition of a noncommutative structure that is
nontrivial only in a vicinity of the diagonal $\Delta_M$ in $M \times
M$. For technical reasons, we will use \emph{formal star products},
i.e. we equip the algebra of functions $C^\infty(M \times M)$ (which
serve as the theory's observables) with a noncommutative product given
by a star product $\tstar$.

Recall that a star product \cite{bayen.et.al:1978a} on a manifold $N$
is a formal $\mathbb{C}[[\lambda]]$-bilinear associative deformation
$\star$ of the algebra of smooth functions $C^\infty(N)$ written as
\begin{equation}
    \label{eq:fstarg}
    f \star g = \sum_{r=0}^\infty \lambda^r C_r (f, g),
\end{equation}
where $f, g \in C^\infty(N)[[\lambda]]$, the $C_r$ are bidifferential
operators such that $\star$ is associative and $1 \star f = f = f
\star 1$ for all $f$.  Additionally, we want $\star$ to be a
deformation of the ordinary product of functions in the sense that
$C_0(f, g) = fg$ is the undeformed commutative product. As a
consequence of associativity, $\{f, g\} = \frac{1}{\I} \left(C_1(f, g)
    - C_1(g, f)\right)$ defines a Poisson bracket and thus a Poisson
bivector field $\theta \in \Gamma^\infty(\Anti^2 TN)$ by $\{f, g\} =
\SP{\theta, \D f \otimes \D g}$. The Jacobi identity for $\{\cdot,
\cdot\}$ is equivalent to $\Schouten{\theta, \theta} = 0$ where
$\Schouten{\cdot, \cdot}$ denotes the Schouten-Nijenhuis bracket of
multivector fields. Recent reviews on deformation quantization can be
found in \cite{gutt:2000a, dito.sternheimer:2002a}. For an elementary
introduction see \cite{waldmann:2006a:script}.

In our case, we consider a star product $\tstar$ on $M \times M$ whose
first order term yields a {Poisson bivector} $\tilde{\theta} \in
\Gamma^\infty(T(M \times M))$ on $M \times M$ (conversely, any such
Poisson bivector can be `quantized' into a star product $\tstar$).
The deformation parameter $\lambda$ is a formal parameter, but may be
thought of as a Planck area in our context. Our crucial requirement
now is that the support of $\tilde{\theta}$ be close to the diagonal
$\Delta_M\subset M\times M$, such that the product of functions
differs from the pointwise one only at small distances.

\begin{figure}
    \begin{center}
        \includegraphics[width=16cm]{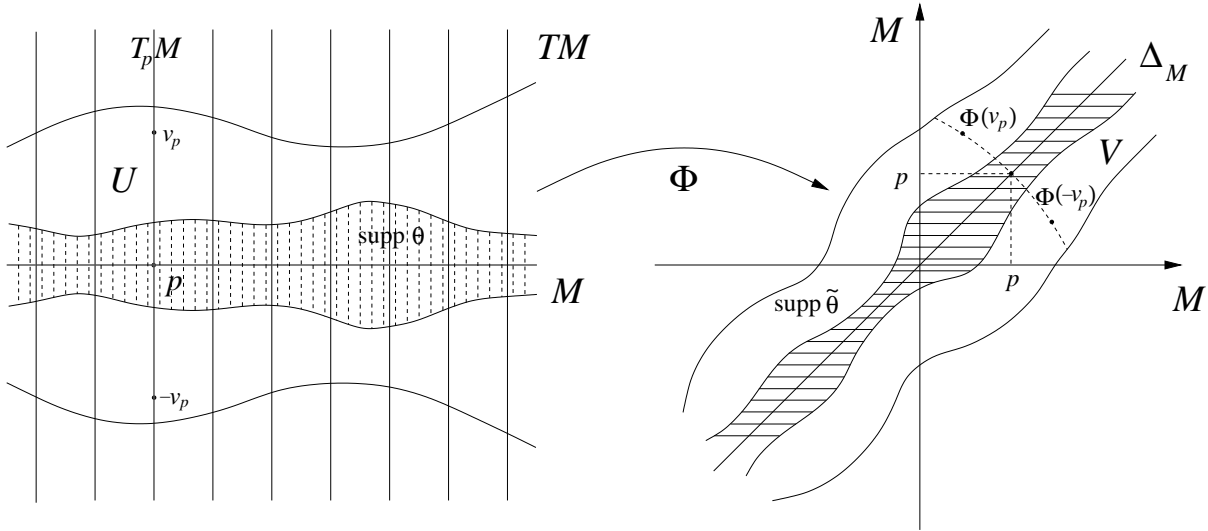}
        \caption{
          \label{fig:SuppThetaTM}
          The exponential map $\Phi$ transports $\theta$ from $TM$ to
          $\tilde{\theta}$ on $M \times M$.  }
    \end{center}
\end{figure}

First, we require that $\supp \tilde{\theta} \subseteq \mathcal{V}
\subseteq M \times M$. This allows to pull-back $\tilde{\theta}$ via
the diffeomorphism $\Phi$ in order to obtain a Poisson bivector
$\theta = \Phi^* \tilde{\theta} \in \Gamma^\infty(\Anti^2
T\mathcal{U})$. Since $\supp \tilde{\theta} \subseteq \mathcal{V}
\subseteq M \times M$, the Poisson bivector $\theta$ extends to a
globally defined Poisson bivector on $TM$ with support $\supp \theta
\subseteq \mathcal{U}$, see also Figure~\ref{fig:SuppThetaTM}.
Secondly, we require that $\supp \theta \cap T_pM$ is \emph{compact}
in $T_pM$. This expresses in a purely topological manner that the
support of $\theta$ is `small': instead of using the (Lorentz) metric
explicitly, we simply fix the 'range of noncommutativity' (given by
the support of $\theta$, or $\tilde \theta$ respectively) to remain
finite.

\begin{figure}
    \begin{center}
        \includegraphics[width=8cm]{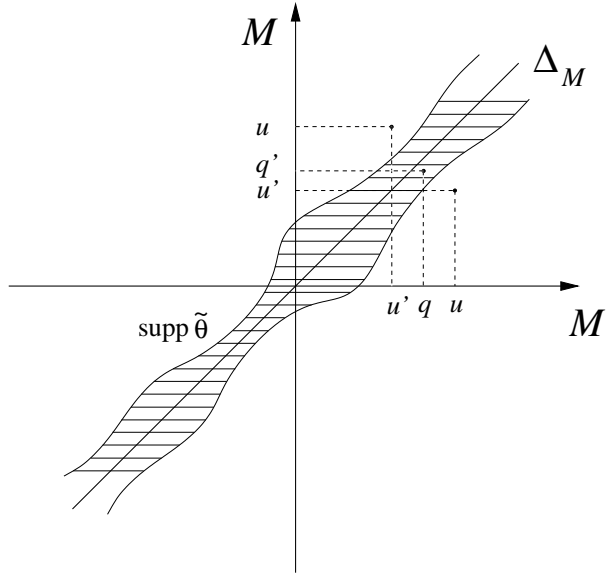}
        \caption{
          \label{fig:SuppTheta}
          Some pairs of points are within the range of
          noncommutativity, some are not. Here the support of
          $\tilde{\theta}$ is symmetric under the canonical flip
          $\tau$.  }
    \end{center}
\end{figure}

The last requirement is to admit only such Poisson bivectors
$\tilde{\theta}$ as are invariant under $\tau_{M \times M}$.
By~\eqref{eq:PhitauInvariant}, this ensures invariance of $\theta$
under $\tau_{TM}$, such that a pair of points $(q,q^\prime)$ is within
the range of noncommutativity if and only if $(q^\prime,q)$ is, see
Figure~\ref{fig:SuppTheta}.

Starting from such a Poisson bivector $\tilde{\theta}$ as first order
term in the star product, the construction methods for star products
like those in \cite{kontsevich:2003a, cattaneo.felder.tomassini:2002b,
  dolgushev:2005a, fedosov:1996a} will yield star products $\tstar$,
whose higher order bidifferential operators $\tilde{C}_r$ for $r \ge
1$ still have support contained in $\supp \tilde{\theta}$. In
principle, there exist more general star products not obeying this
support condition but we shall only use star product with $\supp
\tilde{C}_r \subseteq \supp \tilde{\theta}$.  Clearly, such star
products will reduce to the ordinary pointwise product of functions
outside the support of $\tilde{\theta}$,
\begin{equation}
    \label{eq:ftstargOutsideV}
    (f \tstar g) \Big|_{M \times M \setminus \supp \tilde{\theta}} = 
    (fg) \Big|_{M \times M \setminus \supp \tilde{\theta}}
\end{equation}
for all $f, g \in C^\infty(M \times M)[[\lambda]]$.  Finally, one can
easily arrange that the symmetry $\tau_{M \times M}$ remains a
symmetry for $\tstar$, i.e. that
\begin{equation}
    \label{eq:Reflectiontstar}
    \tau_{M \times M}^* (f \tstar g) = 
    (\tau_{M \times M}^* f) \tstar
    (\tau_{M \times M}^* g)
\end{equation}
for all $f, g \in C^\infty(M \times M)[[\lambda]]$. In the following
we shall always assume that $\tstar$ meets all these requirements. In
fact, we shall discuss even more particular star products and give a
concrete construction for them later.

Thanks to the support properties of $\tstar$ we can pull back each
bidifferential operator $\tilde{C}_r$ to a bidifferential operator
$C_r$ on $\mathcal{U} \subseteq TM$ via the diffeomorphism $\Phi:
\mathcal{U} \longrightarrow \mathcal{V}$. Then these bidifferential
operators yield a star product $\star$ on $\mathcal{U}$ with first
order term corresponding to $\theta$. Thanks to $\supp C_r \subseteq
\supp \theta \subseteq \mathcal{U}$ for all $r \ge 1$, the star
product $\star$ extends to $TM$. Conversely, any star product $\star$
for $\theta$ with the property $\supp C_r \subseteq \supp \theta
\subseteq \mathcal{U}$ for all $r \ge 1$ can be pushed forward via
$\Phi$ to give a star product $\tstar$ on $\mathcal{V}$ which extends
to $M \times M$. Thus both points of view are entirely equivalent as
long as we impose the support conditions. In particular, for $f, g \in
C^\infty(M \times M)[[\lambda]]$ with $\supp f, \supp g \subseteq
\mathcal{V}$ we have
\begin{equation}
    \label{eq:PhiIsomorphismus}
    \Phi^* (f \tstar g) = \Phi^* f \star \Phi^* g.
\end{equation}
The property \eqref{eq:Reflectiontstar} translates into the symmetry
\begin{equation}
    \label{eq:TauTMSymmetry}
    \tau_{TM}^*(f \star g) = \tau_{TM}^* f \star \tau_{TM}^* g
\end{equation}
for all $f, g \in C^\infty(TM)[[\lambda]]$. In the following we shall
use both descriptions and pass from one to the other freely.

%
% Vertical star products
%

\subsection{Vertical star products}
\label{subsec:VerticalStarProducts}

Up to now, the bivectors $\theta$ and $\tilde{\theta}$ as well as the
corresponding star products $\star$ and $\tstar$, respectively, can
still be very general as the support conditions alone are not very
restrictive.

\begin{figure}
    \begin{center}
        \includegraphics[width=8cm]{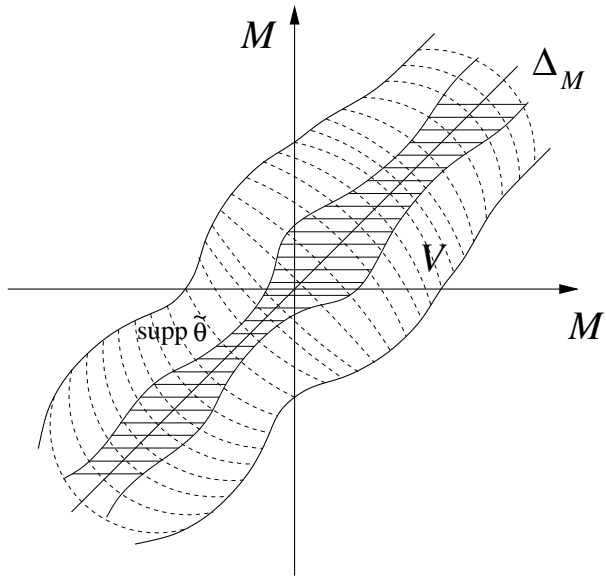}
        \caption{
          \label{fig:RelKoord}
          The dashed lines indicate the locally defined geodesic
          `relative coordinates' transversal to the diagonal. They are
          only defined within $\mathcal{V}$ using the local
          diffeomorphism $\Phi$.  }
    \end{center}
\end{figure}

We now impose one further condition which implements the idea that it
is only the distance between two points that determines whether
noncommutativity is present, while their absolute position in
space-time should not matter (though their absolute position may
influence the specific form of noncommutativity via the bivector's
parametric dependence on $p$). There being {no intrinsic} coordinates
on $M \times M$ transversal to $\Delta_M$, we first use the map $\Phi$
to define \emph{geodesic relative coordinates} near $\Delta_M$, see
also Figure~\ref{fig:RelKoord}. We now ask $\tstar$ to meet the
following additional property: if after restriction to the open subset
$\mathcal{V}$, a function $f \in C^\infty(M \times M)[[\lambda]]$ is
\emph{constant} with respect to the relative coordinates, then for any
other function $g \in C^\infty(M \times M)[[\lambda]]$ we require
\begin{equation}
    \label{eq:fConstanttstargfg}
    f \tstar g = fg = g \tstar f.
\end{equation}
Note that \eqref{eq:fConstanttstargfg} is trivially fulfilled outside
of $\supp \tilde{\theta}$ by \eqref{eq:ftstargOutsideV}. In more
physical terms, observables not sensitive to the relative coordinates
should behave entirely classical, i.e. commutative.
This makes our idea more precise, that $M \cong \Delta_M \subseteq M
\times M$ should remain \emph{commutative} as we have argued in the
introduction: Indeed, functions $f \in C^\infty(M)[[\lambda]]$ can be
prolongated at least locally on $\mathcal{V}$ from $M \cong \Delta_M$
to functions on $M \times M$ by defining them to be constant along the
relative coordinates. Thus a non-trivial star product between such
functions would result in a non-trivial multiplication law for
functions on $M$.

Again, we can translate \eqref{eq:fConstanttstargfg} back to an
equivalent statement regarding $\star$ on $TM$. Here, if $f \in
C^\infty(TM)[[\lambda]]$ is \emph{constant} along the fibers, i.e. of
the form $f = \pi^* u$ with some $u \in C^\infty(M)[[\lambda]]$, then
for any other $g \in C^\infty(TM)[[\lambda]]$ the star product becomes
trivial,
\begin{equation}
    \label{eq:piuStargpiug}
     g \star \pi^* u = g \pi^*u = \pi^* u \star g.
\end{equation}
A more direct characterization of star products with this additional
property is provided by the following theorem:
\begin{theorem}
    \label{theorem:ConditionMeansVertical}
    A star product $\star$ on $TM$ satisfies \eqref{eq:piuStargpiug}
    if and only if $\star$ is vertical, i.e. each bidifferential
    operator $C_r$ differentiates only in vertical directions.
\end{theorem}
Note that for the non-trivial direction of this statement we have to
use the associativity of $\star$. In the
Appendix~\ref{sec:VerticalFormality} we have collected informations on
vertical Poisson structures and vertical star products including their
existence and classification by means of a vertical formality theorem.
As Theorem~\ref{theorem:ExistenceClassification} provides us a
functorial construction of $\star$ out of a given $\theta$
guaranteeing all our requirements, the reader not interested in the
technical details may safely proceed from here on. We summarize the
data and requirements of our model in Table~\ref{tab:Summary}.
\begin{table}[ht]
    \centering
    \begin{tabular}{|c|c|}
        \hline
        Semiclassical $\theta \in \Gamma^\infty(\Anti^2 T(TM))$ 
        & Formal deformation $\star = \sum_{r=0}^\infty \lambda^r C_r$ \\
        \hline
        $\Schouten{\theta, \theta} = 0$ (Jacobi identity)
        & $\star$ associative formal star product\\
        $\theta$ vertical               
        & $\star$ vertical \\
        $\supp \theta \subseteq \mathcal{U}$
        & $\supp C_r \subseteq \supp{C_1} \subseteq \mathcal{U}$ for all
        $r \ge 1$. \\
        $\tau_{TM}^* \theta = \theta$
        & $\tau_{TM}^*$ is automorphism of $\star$ \\
        $\supp \theta \cap T_pM$ compact for all $p \in M$
        & $\supp C_r \cap T_pM$ compact for all $p \in M$ \\
        \hline
    \end{tabular}
    \caption{Summary of the model}
    \label{tab:Summary}
\end{table}

To make ourselves familiar with vertical star products, let us now
give some local formulas. Let $(x^1, \ldots, x^n)$ be local
coordinates on $U \subseteq M$ and denote by $(q^1 = x^1 \circ \pi,
\ldots, q^n = x^n \circ \pi, v^1, \ldots, v^n)$ the induced
coordinates on $TU \subseteq TM$. Here, as usual, a tangent vector
$v_p \in T_pM$ is written as $v_p = v^i(v_p) \frac{\partial}{\partial
  x^i}$, thus specifying the \emph{linear} coordinates $(v^1, \ldots,
v^n)$ on the fibers. Now, a bivector $\theta \in \Gamma^\infty(\Anti^2
T(TM))$ is vertical if and only if locally
\begin{equation}
    \label{eq:ThetaVertical}
    \theta\Big|_{TU} = \frac{1}{2} \theta^{ij}
    \frac{\partial}{\partial v^i} \wedge 
    \frac{\partial}{\partial v^j}\ ,
\end{equation}
where $\theta^{ij} \in C^\infty(TU)$ are local coefficient functions
depending on all variables, $q$'s as well as $v$'s.  The condition
$\tau_{TM}^* \theta = \theta$ means that the functions $\theta^{ij}$
must be even functions of the $v$'s. The support condition is
equivalent to $\theta^{ij}(q, \cdot)$ having compact support with
respect to the $v$-coordinates for fixed $q^1, \ldots, q^n$.  It
follows directly from the definition that the star product $\star$ is
vertical if and only if locally the bidifferential operators $C_r$ are
of the form
\begin{equation}
    \label{eq:VerticalCr}
    C_r(f, g) \Big|_{TU} =
    \sum_{I, J} C_r^{IJ} 
    \frac{\partial^{|I|} f}{\partial v^I}
    \frac{\partial^{|J|} g}{\partial v^J},
\end{equation}
with multi-indices $I$ and $J$, and where the local coefficient
functions $C_r^{IJ} \in C^\infty(TU)$ may again depend on $q$'s as
well as $v$'s. The important point is that both functions are
differentiated \emph{only} in direction of the fiber variables. Of
course, the $C_r^{IJ}$ are subject to further conditions arising from
the associativity of $\star$.

Local expressions for $\tilde{\theta}$ and $\tstar$ are more
complicated as they require knowledge of the explicit form of the
exponential map. This is only in very limited cases accessible whence
we shall mainly work with $TM$ instead of $M \times M$. Note however,
that for the interpretation of functions $f \in C^\infty(TM)$ as
observables one should rather consider their counterparts on $M \times
M$.

%%%%%%%%%%%%%%%%%%%%%%%%%%%%%%%%%%%%%%%%%%%%%%%%%%%
% Further properties of vertical star products
%

\section{Further properties of vertical star products}
\label{sec:PropertiesVertical}

Consider now a vertical Poisson structure $\theta$ and a
corresponding star product $\star$ on $TM$ obeying the support
conditions as well as the reflection symmetry $\tau_{TM}^* \theta =
\theta$ and \eqref{eq:TauTMSymmetry}, respectively. Let
\begin{equation}
    \label{eq:iotapTpMTM}
    \iota_p: T_pM \longrightarrow TM
\end{equation}
denote the embedding of the tangent space at $p \in M$ into the
tangent bundle. As discussed in
Appendix~\ref{subsec:VerticalMultivectors} and
Theorem~\ref{theorem:ExistenceClassification} we can restrict $\theta$
and $\star$ to a Poisson structure $\theta_p$ with corresponding
Poisson bracket $\{\cdot, \cdot\}_p$ and a star product $\starp$ on
$T_pM$, 
\begin{equation}
    \label{eq:RestrictTheta}
    \iota_p^*\left(\{f, g\} \right)
    =
    \left\{ \iota_p^* f, \iota_p^* g \right\}_p
\end{equation}
and
\begin{equation}
    \label{eq:RestrictStarTpM}
    \iota_p^*\left(f \star g \right)
    =
    \iota_p^* f \starp \iota_p^*g
\end{equation}
for all $f, g \in C^\infty(TM)[[\lambda]]$. Here, it is important that
the Poisson structure $\theta$ as well as the star product $\star$ are
vertical, i.e. that all derivatives are only in the direction of the
fibers. By construction, both structures are non-trivial only on
$\mathcal{U}_p$.

In a last step we can push forward both $\theta_p$ and $\starp$ to $M$
via the exponential map $\exp_p$. Since $\supp \theta_p \subseteq
\mathcal{U}_p$ and $\exp_p$ is a diffeomorphism $\exp_p: \mathcal{U}_p
\longrightarrow \exp_p(\mathcal{U}_p) = \mathcal{V}_p$ by our choice
of $\mathcal{U}$, this is well-defined and yields a Poisson bivector
$\tilde{\theta}_p \in \Gamma^\infty(\Anti^2 T\mathcal{V}_p)$.  Again,
the support conditions enable us to extend $\tilde{\theta}_p$ to all
of $M$, whence we obtain a Poisson bivector $\tilde{\theta}_p \in
\Gamma^\infty(\Anti^2 TM)$. Analogously, we obtain a star product
$\tstarp$ on $M$ which quantizes $\tilde{\theta}_p$. Now, by the very
construction of $\tilde{\theta}_p$ and $\tstarp$, for $f, g \in
C^\infty(M)[[\lambda]]$ with $\supp f, \supp g \subseteq
\mathcal{V}_p$, {\sc }
\begin{equation}
    \label{eq:expptildethetap}
    \left\{\exp_p^* f, \exp_p^*g \right\}_{\theta_p}
    =
    \exp_p^* \{f, g\}_{\tilde{\theta}_p}
\end{equation}
and
\begin{equation}
    \label{eq:exppstarp}
    \exp_p^* f \starp \exp_p^* g 
    =
    \exp_p^* \left(f \tstarp g\right).
\end{equation}
This way, every point $p \in M$ obtains its own star product $\tstarp$
being non-trivial only in a neighborhood of the point $p$.

Consider, for example, the linear fiber coordinates $v^i \in
C^\infty(T_pM)[[\lambda]]$, viewed as functions on the tangent space.
Then $v^i \star v^j = v^iv^j + \lambda C_1(v^i, v^j) + \cdots$ whence
$[v^i, v^j]_\star = \I\lambda \theta_p^{ij} + \cdots$. In fact, the
star product $\star$ can be chosen in such a way that for the linear
coordinates $v^i$ the commutator has only the first order terms in
$\lambda$. In any case, note that $\theta^{ij}_p$ is \emph{not} a
constant but a function on $T_pM$ with compact support in
$\mathcal{U}_p$. On $M$, this yields
\begin{equation}
    \label{eq:LocalCommutatorsOnM}
    [x^i,x^j]_{\tstarp} 
    = \I\lambda \tilde\theta_p^{ij} + \cdots
\end{equation}
for the geodesic normal coordinates $x^i$, viewed as local functions
on $M$. Here, the coefficient functions $\tilde\theta_p^{ij}$ are zero
outside a neighborhood of $p$ contained in $\mathcal{V}_p$.

\begin{figure}
    \begin{center}
        \includegraphics[width=14cm]{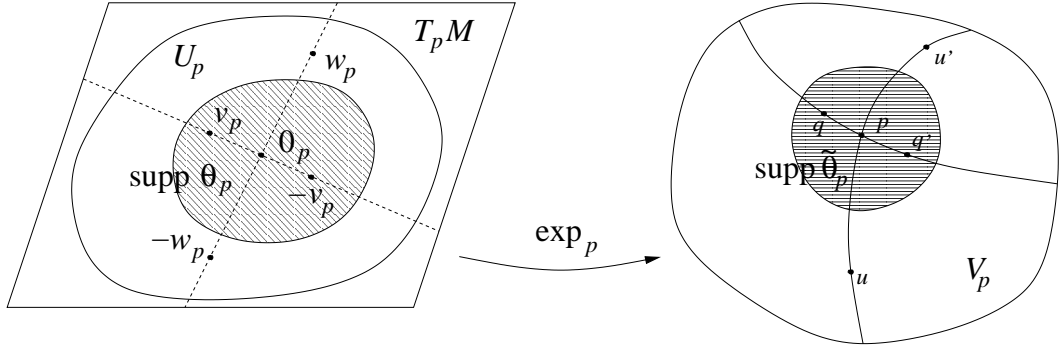}
        \caption{
          \label{fig:NaheBeip}
          The points $q$ and $q'$ are in within the noncommutative
          bubble around their geodesic midpoint $p$, the points $u$
          and $u'$ are still outside.
        }
    \end{center}
\end{figure}

\begin{remark}
    \label{remark:ManyStarProductsOnM}
    Let us emphasize now clearly the interpretation of the star
    products $\tstarp$ compared to the usual star products on the
    space-time manifolds $M$ as mentioned in the introduction. The
    main difference is that we now have a whole family of star
    products $\{\tstarp\}_{p \in M}$ instead of just one. The
    interpretation of the algebra $(C^\infty(M)[[\lambda]], \tstarp)$
    comes from the global picture $(C^\infty(M \times M)[[\lambda]],
    \tstar)$. We are still discussing observables of \emph{two}
    points, i.e. functions $f$ on $M \times M$ and \emph{not}
    functions on $M$.  However, we may be interested in \emph{states},
    like the $\delta$-functionals $\delta_{q, q'}$ on $M \times M$ and
    their quantum analogs, say for $(q, q') \in \mathcal{V}$ to make
    things non-trivial. Then, thanks to verticality, all we have to
    know to evaluate observables in such a state are their
    \emph{restrictions} to the dashed lines in
    Figure~\ref{fig:RelKoord} through their \emph{geodesic midpoint}
    $p$, see Figure~\ref{fig:NaheBeip}. Then we can equivalently work
    with $(C^\infty(M)[[\lambda]], \tstarp)$ and the restriction of
    the observable $f$ to such a dashed line eventually yields a
    corresponding function on $M$ via $\iota_p^*$ and $(\exp_p)_*$.
    Clearly, we have to make precise what notion of states we are
    going to use.
\end{remark}

%
% Hermitean vertical star products and their states
%

\subsection{Hermitean vertical star products and their states}
\label{subsec:HermiteanStates}

The observables of our theory are functions $f \in C^\infty(M \times
M)[[\lambda]]$. We will now specify states and the corresponding
expectation values for these observables, employing the usual
techniques of deformation quantization, which are a straightforward
analogue of the well-known approaches in $C^*$-algebra or
$O^*$-algebra theory, see e.g.\cite{waldmann:2005b} for a review.

First we make the additional assumption that $\tstar$ and hence also
$\star$, $\starp$ and $\tstarp$ are \emph{Hermitean}, i.e. we require
\begin{equation}
    \label{eq:HermiteanStarProduct}
    \cc{f \tstar g} = \cc{g} \tstar \cc{f}
\end{equation}
for all $f, g \in C^\infty(M \times M)[[\lambda]]$, where
$\cc{\lambda} = \lambda$ is treated as a \emph{real} quantity.  If we
construct $\tstar$ out of $\tilde{\theta}$ or, equivalently, $\star$
out of $\theta$ using a formality as in
Theorem~\ref{theorem:ExistenceClassification}, then the reality
$\theta = \cc{\theta}$ of the Poisson structure implies that the
corresponding star product is Hermitean. We can therefore safely
assume \eqref{eq:HermiteanStarProduct} for $\tstar$, $\star$, $\starp$
and $\tstarp$ in the following.

The complex conjugation now being an involution, we define states as
\emph{positive $\mathbb{C}[[\lambda]]$-linear functionals}
\begin{equation}
    \label{eq:PositiveFunctional}
    \omega: C^\infty(M \times M)[[\lambda]]
    \longrightarrow \mathbb{C}[[\lambda]]
    \quad
    \textrm{with}
    \quad
    \omega(\cc{f} \star f) \ge 0,
\end{equation}
where the positivity is understood in the sense of formal power series
(a real formal power series $a = \sum_{r=r_0}^\infty \lambda^r a_r \in
\mathbb{R}[[\lambda]]$ is positive, if the lowest non-vanishing
coefficient $a_{r_0}$ is positive, $a_{r_0} > 0$). In addition to
\eqref{eq:PositiveFunctional} we require that states be
\emph{normalized}, $\omega(1) = 1$.  The
$\mathbb{C}[[\lambda]]$-linearity implies that $\omega$ is of the form
\begin{equation}
    \label{eq:omegaSumomegar}
    \omega = \sum_{r=0}^\infty \lambda^r \omega_r
    \quad
    \textrm{with $\mathbb{C}$-linear maps}
    \quad
    \omega_r: C^\infty(M\times M) \longrightarrow \mathbb{C}.
\end{equation}
In particular, $\omega_0$ turns out to be a positive
$\mathbb{C}$-linear functional of the commutative $^*$-algebra
$C^\infty(M \times M)$, i.e. $\omega_0(\cc{f}f) \ge 0$ for all $f \in
C^\infty(M \times M)$. It follows that $\omega_0$ is the integration
with respect to a compactly supported positive Borel measure on $M
\times M$.

Conversely, and this is the important point here, one can show that
\emph{any} classical $\omega_0$ can be deformed into a functional
$\omega$ which is a state with respect to $\tstar$, see
\cite{bursztyn.waldmann:2005a}.  Note that the `quantum corrections'
$\omega_r$ to $\omega_0$, which, in general, are necessary to ensure
positivity, are \emph{by no means unique}: there are many quantum
states $\omega$ with the same classical limit $\omega_0$. Although it
is generally very difficult to find the corrections explicitly, one
can show that they can always be chosen to be of the form $\omega_r =
\omega_0 \circ S_r$ with a differential operator $S_r$. In such a
case, the support of $\omega$ coincides with that of $\omega_0$.

Due to the positivity of $\omega$, we may interpret $\omega(f)$ as the
expectation value of the observable $f$ in the state $\omega$. Now,
$\omega$ still satisfies a Cauchy-Schwarz inequality (in the sense of
formal power series), allowing us to write down uncertainty relations.
As usual, we define the variance of an observable $f$ in the state
$\omega$ by
\begin{equation}
    \label{eq:Variance}
    \Var_\omega(f) 
    = 
    \omega\left( \cc{(f - \omega(f))} \tstar (f - \omega(f))\right) 
    \ge 0.
\end{equation}
Then for two Hermitean elements $f = \cc{f}$ and $g = \cc{g}$, i.e.
observables in the stricter sense, we find the usual uncertainty
relation
\begin{equation}
    \label{eq:Uncertainty}
    4 \Var_\omega(f) \Var_\omega(g) \ge
    \omega\left([f, g]_{\tstar}\right)
    \cc{\omega\left([f, g]_{\tstar}\right)}\ ,
\end{equation}
where the Cauchy-Schwarz inequality for $\omega$ has been used.  As
usual, these inequalities justify the identification of positive
functionals with states.
\begin{remark}
    \label{remark:FormalPositivity}
    We should note that the notion of positivity we are using is on
    one hand the only reasonable from an algebraic point of view: it
    is the unique one which makes $\mathbb{R}[[\lambda]]$ an ordered
    ring such that $\lambda > 0$. On the other hand, there is a more
    concrete motivation coming from asymptotics: if we think of our
    formal star product $\star$ as being the asymptotic expansion of
    some convergent product, say in a $C^*$-algebraic approach, then
    one can also asymptotically expand positive linear functionals
    which yield precisely the ones we are studying. In this sense, the
    notion of positivity we are using is the best we can have. See
    also \cite{waldmann:2005b} for a more detailed discussion of
    states in deformation quantization.
\end{remark}
\begin{remark}
    \label{remark:ObservablesUnchanged}
    It will be important for the physical interpretation to note how
    the noncommutative structure has entered here: $\omega$ is a
    deformed classical state $\omega_0$ whose quantum corrections
    depend on $\tstar$, such that the expectation values of the
    observable $f$ is changed when we pass from classical to
    noncommutative space-time. It is the main feature of the
    deformation approach that the observable $f$ itself remains
    \emph{unchanged}: it is still the same function with the same
    physical interpretation as observable. We only changed the product
    structure and hence the states.
\end{remark}

One effect of noncommutativity is that the variances of observables
will in general be strictly larger than the classical ones. In
particular, the $\delta$-functionals $\delta_{(q, q')}$ for $(q, q')
\in M \times M$ are no longer positive with respect to $\tstar$, but
require quantum corrections, and we will always find observables such
that the variances in these deformed $\delta$-functionals are strictly
positive (while the classical ones are of course 0). We shall come
back to explicit examples in
Section~\ref{sec:NoncommutativeMinkowski}.

Let us now discuss why our model meets the physical requirements which
we have argued for. We consider now a classical state $\omega_0$, i.e.
a positive Borel measure on $M \times M$ whose (compact) support
$\supp \omega_0$ is far away from the diagonal $\Delta_M$, and in
particular, $\supp \omega_0 \cap \supp \tilde{\theta} = \emptyset$. It
immediately follows that $\omega_0$ is a state with respect to
$\tstar$, since $\tstar$ is non-trivial only in $\supp
\tilde{\theta}$. Indeed, we have
\begin{equation}
    \label{eq:omegaNullfstargfg}
    \omega_0(f \tstar g) = \omega_0(fg)
\end{equation}
for all $f, g \in C^\infty(M \times M)[[\lambda]]$ in this case.
Therefore, all variances and covariances of $f$ and $g$ with respect
to $\omega_0$ are the classical ones. In particular, only the
classical variances appear in \eqref{eq:Uncertainty} and the right
hand side is zero, although $[f, g]_{\tstar}$ may be different from
zero.  This shows that if we evaluate observables $f \in C^\infty(M
\times M)[[\lambda]]$ far away from the diagonal, no noncommutative
behavior can be seen. The noncommutativity only appears close to the
diagonal as is expected from the support conditions on
$\tilde{\theta}$ and $\tstar$. This is precisely the behavior we
wanted. At large distances our locally noncommutative space-time
behaves entirely classically.

%%%%%%%%%%%%%%%%%%%%%%%%%%%%%%%%%%%%%%%%%%%%%%%%%%%

% Distance measurements and the causal structure

\subsection{Distance measurements and the causal structure}
\label{subsec:distMeas}

Let us now reconsider the interpretation of our noncommutative
structure from the point of view of distance measurements.  Since the
concept of `distance' is of course misleading in a pseudo-Riemannian
context, we shall not measure a distance function, but measure the
metric directly. It turns out that this can be done most natural in
our framework.

As a motivation one may think of a Riemannian situation where the
metric distance $d(q, q')$ between two points $q$ and $q'$ is defined
as the infimum over the lengths of all paths joining the two points.
In general, this is a highly non-trivial quantity. However, if the
points are close enough then one finds a unique shortest geodesic
joining them, whose length realizes $d(q, q')$. In fact, if $(q, q')
\in \mathcal{V}$ then this is the case and the geodesic is precisely
the one starting from the geodesic midpoint $p$ in opposite directions
where $- v_p = \exp_p^{-1}(q)$ and $v_p = \exp_p^{-1}(q')$. In this
case, the distance is given by $d(q, q') = 2 \sqrt{g_p(v_p, v_p)}$. In
particular, the \emph{square} of the distance function is the
\emph{smooth} function $d^2(q, q') = 4 g_p (v_p, v_p)$. In general,
the distance function is only smooth close to the diagonal.

In the general situation we shall therefore use the function $d^2 \in
C^\infty(TM)$ defined by
\begin{equation}
    \label{eq:dsquare}
    d^2(v_p) = g_p(v_p, v_p)
\end{equation}
as a good replacement for the geodesic distance function. It is a
quadratic function on $TM$ which is everywhere smooth and in the
neighborhood $\mathcal{U}$ it is indeed the square of the distance
function in the Riemannian case. Since we are only interested in the
behavior close to the diagonal $\Delta_M$ as the noncommutativity is
only present here, this will be a perfect observable to measure the
metric.

The quantum effects will now come into the game in the expectation
values of this observable $d^2$ if we evaluate it in some state. In
particular, we are interested in those states which are as close as
possible to the $\delta$-functionals at some point $v_p$. Thanks to
our verticality condition we can consider even the restricted
situation, i.e. the observable $\iota_p^* d^2 \in
C^\infty(T_pM)[[\lambda]]$. Then we need a deformation of the
$\delta$-functional
\begin{equation}
    \label{eq:PositiveDelta}
    \delta_{v_p}^{(\starp)} 
    = \delta_{v_p} \circ S_{v_p},
    \quad
    \textrm{where}
    \quad
    S_{v_p} = \id + \sum_{r=1}^\infty \lambda^r S_{v_p}^{(r)},
\end{equation}
into a positive functional for $\starp$. Then the quantum distance
square between $q$ and $q'$ is now (up to the factor $4$) the
evaluation
\begin{equation}
    \label{eq:ExpectationValueDistance}
    \delta_{v_p}^{(\starp)} (d^2) 
    = d^2(v_p) 
    + \lambda \left(S^{(1)}_{v_p} d^2\right)(v_p) + \cdots,
\end{equation}
which is clearly a deformation of the classical distance square.
Moreover, in general we obtain a non-trivial variance of this
measurement according to \eqref{eq:Variance} since on one hand $d^2
\starp d^2$ is not just the pointwise product and on the other hand
the correction terms $S^{(r)}_{v_p}$ contribute as well. This way we
arrive at the observation that the geometry indeed becomes
\emph{fuzzy}. Note however, that the choice for a deformation of
$\delta_{v_p}$ is not unique at all.

Let us also remark already at this point that in the Lorentz situation
the sign of the classical evaluation $d^2(v_p)$ determines whether
$v_p$ is a space-like, light-like or time-like vector. In our case,
this characterization needs no longer to be valid, in particular, the
light-like vectors with $d^2(v_p) = 0$ might get correction terms from
the deformed $\delta$-functional which makes them space-like or
time-like. Note however, that this again depends of course on our
choice of the deformation $\delta_{v_p}^{(\starp)}$: this simply
reflects again that there are no `classical' points any more in a
truly noncommutative space-time.

%%%%%%%%%%%%%%%%%%%%%%%%%%%%%%%%%%%%%%%%%%%%%%%%%%%%%%%%%%%%
% The Noncommutative Minkowski Space
%

\section{The Noncommutative Minkowski Space}
\label{sec:NoncommutativeMinkowski}

\begin{figure}[htbp]
    \begin{center}
        \includegraphics[width=8cm]{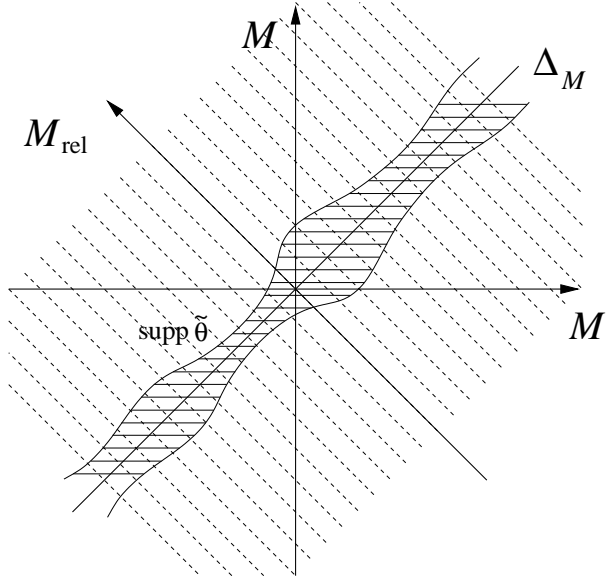}
        \caption{
          \label{fig:MinkowskiI}
          In the case of Minkowski space the map $\Phi$ is a global
          diffeomorphism and introduces the global center of mass and
          relative coordinates on $M \times M$.
        }
    \end{center}
\end{figure}

In order to analyze our construction's properties more explicitly, we
now discuss the case where $M$ is a vector space of dimension $n$ in
more detail. We chose $\nabla$ to be the canonical flat connection.
Then the exponential map at each point is a global diffeomorphism, and
so is $\Phi$.  In fact, the exponential map implements a
diffeomorphism of $T_pM$ to $M$ given by the translation,
\begin{equation}
    \label{eq:ExpIsAddition}
    \exp_p(v)=p+v \in M \qquad \forall\ v\in T_pM,
\end{equation}
and the local formulae of the previous section are now globally
defined. In particular, for any pair of points $(q, q^\prime) \in M
\times M$, the midpoint $p=(q + q^\prime)/2$ and the relative
coordinates $v = (-q + q^\prime)/2$ are now globally defined.  For
this reason, the dashed lines from Figure~\ref{fig:RelKoord} which
denote the geodesic relative coordinates and in the general case are
defined only in some open neighborhood $\mathcal V\subset M\times M$
of the diagonal $\Delta_M$, now become straight lines extending to
infinity, see Figure~\ref{fig:MinkowskiI}. Note also, that up to now
we have only used the connection $\nabla$ but no metric.

Let $f,g \in C^\infty (TM)[[\lambda]]$, then
\begin{equation}
    \label{eq:fStargAgain}
    f\star g = \sum \lambda^r C_r(f,g),
\end{equation}
where
\begin{equation}
    \label{eq:CrLocallyOnMinkowski}
    C_r(f,g) =
    \sum_{I,J} C_r^{IJ}
    \frac{\partial^{|I|} f}{\partial v^I} 
    \frac{\partial^{|J|} g}{\partial v^J},
\end{equation}
with multi indices $I, J \subset \{1,\dots,n\}$, $|I|,|J|\leq r$. Note
that $C_r$ differentiates only in direction of the tangent spaces,
i.e. in $v$-directions.  Thanks to the simple form of the exponential
map \eqref{eq:ExpIsAddition} we find from \eqref{eq:PhiIsomorphismus}
for all $f, g \in C^\infty(M\times M)[[\lambda]]$
\begin{equation}
    f \tstar  g 
    = \sum \lambda^r 
    \sum_{I,J} \tilde{C}_r^{IJ} 
    \prod_{i\in I}
    {\textstyle \left(-\frac{\partial}{\partial q^i}
          +\frac{\partial}{\partial {q^\prime}^i}\right)}
    f
    \prod_{j\in J}
    {\textstyle \left(-\frac{\partial}{\partial q^j}
          +\frac{\partial}{\partial {q^\prime}^j}\right)} 
    g,
\end{equation}
where $\tilde{C}_r^{IJ} = C_r^{IJ}\circ \Phi^{-1}$. Here, we have used
that $\partial_v (\Phi^*f)(p,v)=\big((-\partial_1+\partial_2) f\big)
(\Phi(p,v))$ with $\partial_i$ denoting the derivative with respect to
the $i^{th}$ argument.

%%%%%%%%%%%%%%%%%%%%%%%%%%%%%%%%%%%%%%%%%%%%%%%%%%%%%
% Global constant Poisson structure
%

\subsection{Global constant Poisson structure}

As $\Phi$ is a global diffeomorphism on flat space, there is in
principle no need to restrict the range of noncommutativity, i.e. to
have $\tilde{\theta}$ compactly supported. We will however, of course
choose to still impose such restrictions in order to implement
localized noncommutativity, see Figure~\ref{fig:MinkowskiI}. But for
the time being, in order to compare our approach to more commonly
analyzed settings, we now restrict ourselves to the special case of a
\emph{constant} vertical Poisson structure
\begin{equation}
    \label{eq:ThetaVeryConstant}
    \theta = \frac{1}{2} \theta^{ij} 
    \frac{\partial}{\partial v^i} \wedge \frac{\partial}{\partial v^j}
    \quad
    \textrm{with}
    \quad
    \theta^{ij} = - \theta^{ji} \in \mathbb{R}.
\end{equation}
Then, the star product on $TM$ can be chosen to be the usual
Weyl-Moyal product
\begin{equation}
    \label{eq:Constantstar}
    f \star g 
    = \sum_{r=0}^\infty
    \frac{1}{r!}
    \left(\frac{\I\lambda}{2}\right)^r
    \theta^{i_1 j_1} \cdots \theta^{i_r j_r} 
    \frac{\partial^r f}{\partial v^{i_1} \cdots v^{i_r}}
    \frac{\partial^r g}{\partial v^{j_1} \cdots v^{j_r}}.
\end{equation}
It is obviously invariant under reflections in the sense of
\eqref{eq:TauTMSymmetry}. With such a constant Poisson structure, any
pair of points is within the range of noncommutativity of their
midpoint, since $\supp \tilde \theta_p = T_pM$.

Clearly, the dependence on $p$ is only in the functions $f$ and $g$
and we recognize that \eqref{eq:Constantstar} restricts to the
Weyl-Moyal star product on $T_pM$ with respect to $\theta$, i.e.
\begin{equation}
    \label{eq:WeylMoyalOnTpM}
    f \starp g 
    = \mu \circ \exp 
    \left( 
        -\frac{\I\lambda}{2} \theta^{ij} 
        \partial_{v^i} \otimes \partial_{v^j} 
    \right)
    (f \otimes g)
\end{equation}
for $f, g \in C^\infty(T_pM)[[\lambda]]$ where $\mu(f\otimes g)= fg$
is the usual pointwise product.  Although in fact, the star product is
independent of $p$, we keep the notation $\starp$ in order to
remember that we are considering some fixed $T_pM$.  Likewise, we
find for $f,g \in C^\infty (M\times M)[[\lambda]]$
\begin{equation}
    \label{eq:Constantstarq}
    f \tstar  g
    = 
    \mu \circ \exp
    \left(
        \frac{\I\lambda}{2} \theta^{ij}
        \left(- \frac{\partial}{\partial q^i} +
            \frac{\partial}{\partial {q'}^i}
        \right)
        \otimes
        \left(- \frac{\partial}{\partial q^j} +
            \frac{\partial}{\partial {q'}^j}
        \right)
    \right)
    (f \otimes g).
\end{equation}
Obviously, the differentiation is in the direction of the line through
$q$ and $q^\prime$ (i.e. perpendicular to the diagonal).

As in Section~\ref{subsec:distMeas} we now ask ourselves how close two
points $(q, q^\prime) \in M \times M$ may be to one another. In order
to do so, we consider the situation on $T_pM$, and modify the
$\delta$-distribution such that it is a positive functional with
respect to the star product $\star_p$. As discussed in
Section~\ref{subsec:HermiteanStates}, this deformation is not unique,
but as a natural candidate we use the formal version of the coherent
states of quantum mechanics, see the discussion in
\cite{bursztyn.waldmann:2000a, bursztyn.waldmann:2005a}. For
simplicity we assume that $\theta$ is non-degenerate, i.e. a
symplectic Poisson tensor. Hence in particular, $M$ has to be even
dimensional. Then we consider
\begin{equation}
    \label{eq.deltaPositive}
    \delta^{(\star_p)}
    = \delta\circ \E^{\frac{1}{4} \lambda \Delta_g}
\end{equation}
and likewise for the translates of the $\delta$-distribution
$\delta_w$, $w\in M$. Here, $g^{-1}$ is a positive compatible scalar
product with $\theta$, i.e. there exists a linear complex structure $J
\in \End(T_pM)$, $J^2 = -\id$ with $g^{-1} (v, w) = \omega_\theta(v,
Jw)$ for all $v, w \in T_pM$, where $\omega_\theta$ is the associated
symplectic form to $\theta$. Finally, $\Delta_g$ denotes the usual
Laplacian with respect to $g^{-1}$. In the following we shall mainly
consider the standard symplectic form $\theta$ and chose for $g^{-1}$
the identity matrix with respect to some given choice of Darboux
coordinates on $T_pM$.

Now, for any quadratic form $A \in M(n,\R)$, $f_A (v) = v^tAv$ and for
a Laplacian $\Delta_g$ with respect to some symmetric form $g^{-1}$,
we find
\begin{equation}
    \label{eq:exponfA}
    \E^{\frac{1}{4} \lambda \Delta_g} f_A = f_A + \tsf \lambda 2 \tr(gA) 
\end{equation}
and, after a short calculation,
\begin{equation}
    \label{eq:exponfAfA}
    \E^{\frac{1}{4} \lambda \Delta_g} (f_A\star f_A)
    =
    f_A^2 + \lambda f_A \tr gA + 2\lambda f_{AgA}
    + \frac{\lambda^2}{4} \left(
        2 \tr(A_\theta A_\theta) + (\tr gA)^2 + 2\tr gAgA
    \right),
\end{equation}
where $(A_\theta)^r_j = \theta^{rs}A_{sj}$ and where we have used
$\Delta_g(f_Af_A) = 4 f_A (\tr gA) + 8f_{AgA}$ and $\Delta_g^2
(f_Af_A) = 8 (\tr gA)^2 + 16 \tr gAgA$.

\begin{figure}[htbp]
    \begin{center}
        \includegraphics[width=8cm]{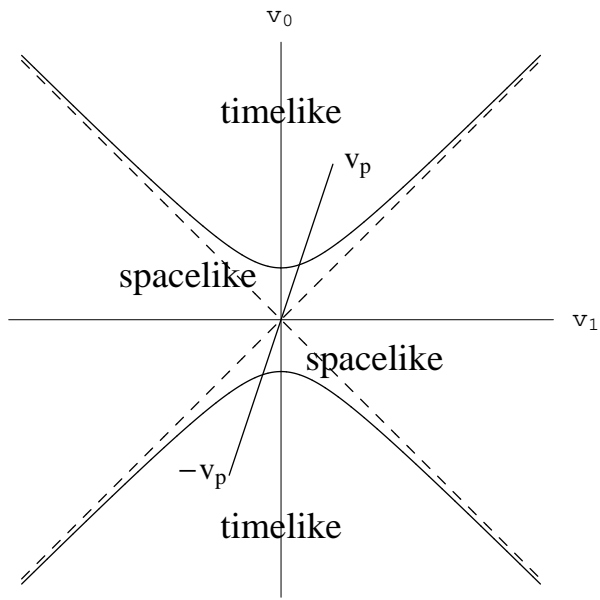}
        \caption{
          \label{fig:defLC}
	  The ordinary and the deformed light-cone in the case of constant
          $\theta$. Two spatial dimensions of $v\in T_pM$ are suppressed.}
    \end{center}
\end{figure}

Following our general discussion in Section~\ref{subsec:distMeas}, let
us now investigate the Lorentz square in 4 dimensions, i.e.  consider
$f_\eta$ where $\eta={\rm diag}\,(+,-,-,-)$ and $n = 4$. By
\eqref{eq:exponfA}, a distance measurement as above yields for the
squared distance,
\begin{equation}
    \label{eq:Deltafeta}
    \delta^{(\star_p)}_{v_p}(f_\eta) 
    =  f_\eta(v_p) + \frac{\lambda}{2} \tr(\eta g)
    =  \eta(v_p, v_p) - \lambda,
\end{equation}
whence all Lorentz squares acquire a {\em negative} offset independent
of $v_p$ in this particular deformation of the classical
$\delta$-functional.  The resulting deformed light cone on $T_pM$,
defined by $\delta^{(\star_p)}_{v}(f_\eta) = 0$, then takes the form
of two hyperbolae, $v_0 = \pm \sqrt{\lambda +\vec v^2}$ (i.e. a ``mass
shell'' of mass $\lambda$), approaching the ordinary light-cone for
distances $\|v\| \gg \lambda$, see Figure~\ref{fig:defLC}.
The time-like vectors are characterized by $v_0^2 > \lambda + \vec
v^2$, and the space-like vectors are those with $v_0^2 < \lambda + \vec
v^2$. To interpret this picture we recall that the point $p$ does not
have meaning in itself, but only as the geodesic midpoint of two other
points $\exp_p(\pm v_p) = p \pm v_p$, such that the points $v_p$ and
$-v_p$ are connected by a time-/space- or light-like line (in the
deformed sense). The gap between future and past time-like lines
around $0_p$ is of the order $\lambda$ and is to be interpreted as
follows: if the two points in $M$, $p \pm v_p$, approach each other,
the causal structure is lost. This is however, by construction in
accordance with our minimal resolvable distance.

Note that the offset would have opposite sign, had we used $\eta =
\mathrm{diag}(-, +, +, +)$, so also in this case, the light-cone is deformed in
the \emph{same} manner as above. For the variance we find
\begin{equation}
    \label{eq:Variance}
    \Var_{\delta^{(\star_p)}_{v_p}}(f_\eta) 
    = 
    f_{\eta^2}(v_p) + 2\lambda^2
\end{equation}
by \eqref{eq:exponfAfA} and the fact that $f_\eta \starp f_\eta =
f_\eta f_\eta$ since $\sum_{r,j}(\eta_\theta)^r_j (\eta_\theta)^j_r =
-\sum (\theta^{rj})^2\,\eta_{jj}\eta_{rr} = 0$. 
\begin{remark}
    \label{remark:Distance}
    We also would like to note that the measurements of the distance
    square are not an artifact of our approach but an intrinsic
    feature of the noncommutative Minkowski space-time with constant
    $\theta$.
\end{remark}

%%%%%%%%%%%%%%%%%%%%%%%%%%%%%%%%%%%%%%%%%%%%%%%%%
% Constant along fiber

\subsection{Non-constant Poisson structures}
\label{subsec:NonConstantPoisson}

In a slightly more general scenario, we might want to employ a Poisson
structure which is constant along each fiber $T_pM$, but varies
depending on the absolute position of the center of mass, $p$. In this
case, the formulae \eqref{eq:Constantstar} and
\eqref{eq:Constantstarq} from the previous discussion remain valid;
the only difference being that $\theta^{ij}$ now explicitly depends on
$p$. This however, seems to be a rather unnatural scenario, since
translation invariance is unnecessarily broken.  In any case, we set
out to construct a noncommutative structure that vanishes in the limit
of large distances. In particular, we require that an appropriate
deformation $\delta^{(\tstar)}_{q,q^\prime}$ of the
$\delta$-Distribution on $C^\infty(M\times M)$, fulfills
(\ref{eq:omegaNullfstargfg}) for $(q,q^\prime)$ far away from the
diagonal $\Delta_M$, i. e. that for $f, g \in C^\infty(M \times
M)[[\lambda]]$, $\delta^{(\tstar)}_{q,q^\prime} (f\tstar g) =
\delta_{q, q^\prime}(fg)$ (or equivalently, on $T_pM$, that for $f, g
\in C^\infty(T_pM)[[\lambda]]$, $\delta^{(\star_p)}_{v_p} (f\star_p g)
= \delta^{(\star_p)}_{v_p}(fg)$ for large $v_p$). Obviously, this is
not true for a nontrivial Poisson structure that is constant along the
fiber $T_pM$ whether it depends on $p$ in a non-trivial way or not.

%\subsection{Compactly supported Poisson structure}

Let us therefore now turn to a scenario which actually exhibits the
features our more general approach allows for and choose $\supp
\theta_p$ to be compact. For concreteness' sake we may think of the
special Poisson structure for which $\theta_p$ is constant on a ball
$\overline{B_r(0_p)}$ around $0_p \in T_pM$ and then decreases quickly
to 0, such that $\supp \theta_p \subset
\overline{B_{r+\epsilon}(0_p)}$ for some $\epsilon$. To implement flip
symmetry (\ref{eq:Reflectiontstar}) we moreover impose that
$\theta_p^{ij}(v_p)=\theta_p^{ij}(-v_p)$ for all $v_p \in T_pM$ (on
$\overline{B_r(0_p)}$ this is of course trivially fulfilled). Such
Poisson structures exist, see
Example~\ref{example:VerticalAndCompactSupp}, and meet all our
requirements. Clearly, we can use the same such $\theta_p$ for all $p$
whence we easily can implement translation invariance. In particular,
the support $\supp \theta \cap T_pM$ is of the same size for all $p
\in M$.

\begin{figure}[htbp]
    \begin{center}
       \includegraphics[width=8cm]{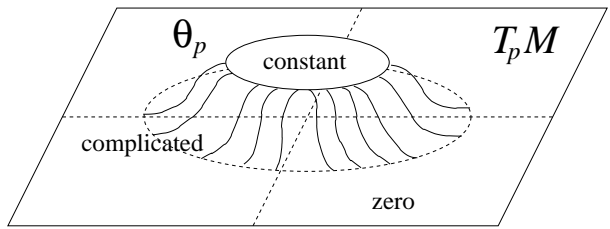}
       \caption{
          \label{fig:SupportTheta}
          Schematic view of a Poisson structure $\theta_p$ being
          constant around $0_p$ with compact support
        }
    \end{center}
\end{figure}

In such a scenario we again consider a distance measurement.  Although
the deformation $\delta_{v_p} \circ \E^{\frac 1 4 \lambda \Delta_g}$
we previously employed, may no longer be a positive functional for
$(C^\infty(T_pM)[[\lambda]],\starp)$ for all classical points $v_p \in
T_pM$, it will be positive for $v_p$ contained in
$\overline{B_{r}(0_p)}$.  Here, $\theta_p$ is constant, and the star
product coincides with the Weyl-Moyal star product as in
\eqref{eq:Constantstar}. If we are interested only in smallest
distances, this set of states is sufficient and evaluating the
function $f_\eta$ in $\delta_{v_p}^{\star_p}$ with $v_p \in B_r(0_p)$,
we gain the same results as in the case of globally constant
$\theta_p$ above. Of course, if we consider states corresponding to
$\delta$-functionals $\delta_{v_p}$ for larger $v_p$, we will have to
consider some other deformation.  Moreover, in the area where
$\theta_p$ drops to zero as a function of $v_p$, the additional
derivatives of $\theta_p$ will contribute significantly to the
distance measurement, see also Figure~\ref{fig:SupportTheta}.
Thus, it will be of major importance to understand the state space of
formal star products better. Hence, a priori the passage from very
small to large distances is not yet very well controlled in such
models. We hope to come back to this question at a later stage.

%%%%%%%%%%%%%%%%%%%%%%%%%%%%%%%%%%%%%%%%%%%%%%%%%%%%%%%%%%%%%%

\subsection{Symmetries}
\label{subsec:Symmetries}

It is natural now to consider the behavior of our construction with
respect to Poincar\'e transformations. In the discussion above we have
already seen that translation invariance can easily be accounted for
by choosing $\theta_p$ and $\starp$ to be independent of the geodesic
midpoint $p$, i.e. to choose the same structure on \emph{all} $T_pM$.
It is furthermore quite simple to implement invariance of the star
product under orthogonal transformations with respect to some
\emph{positive definite} scalar product, simply by asking that
$R^*\theta_p = \theta_p$ for all $R \in O(4)$. Note that
reflections at $0_p$ are already taken care of by implementing the
flip symmetry.

\emph{Lorentz symmetry} on the other hand, cannot be implemented in
such a simple manner. The reason is that on one hand the diagonal
action of the Lorentz group on $M \times M$ induces the usual action
on the relative coordinates in $T_pM$. On the other hand, it is
well-known that there is no Lorentz invariant antisymmetric bivector
on $T_pM$ beside $\theta_p = 0$. Thus we necessarily break Lorentz
invariance already on the semi-classical level of $\theta_p$.

For a \emph{constant} and also translation invariant $\theta_p$ we
might however mimic the approach taken in
\cite{doplicher.fredenhagen.roberts:1995a} and consider along with
some fixed $\theta_p$ the whole orbit under Lorentz transformations.
Thus with this covariant transformation law for $\theta_p$ we would
obtain an action of the Lorentz group on the whole algebra but of
course we have introduced this additional orbit which affects the
classical limit in a non-trivial
way~\cite{doplicher.fredenhagen.roberts:1995a}.

%
% Further Questions and Outlook
%

\section{Further Questions and Outlook}
\label{sec:QuestionOutlook}

Let us conclude with some open questions and further remarks.

%
% Dynamics and field theories
%

\subsection{Dynamics and Field Theories}
\label{subsec:Dynamics}

Up to now we have only set up a kinematic framework for a locally
noncommutative space-time. This is of course not in the least enough
to have a reasonable model for space-time at small distances: We
certainly have to include some sort of dynamics into our description.
Here one should try to proceed in the usual stages.

A first approach would be to investigate the behaviour of point-like
classical or quantum mechanical particles moving in such a locally
noncommutative space-time.  Especially for non-relativistic
considerations, it seems reasonable to consider as a first step only a
locally noncommutative \emph{space} and treat time as an ordinary real
parameter.

In a second step, one can consider classical field theories on our
locally noncommutative space-time. A good starting point is provided
by deformed vector bundles in the sense of \cite{waldmann:2005a,
  waldmann:2001b, bursztyn.waldmann:2000b} concerning the matter part.
It is however, not yet clear how to define interaction terms, though
one might hope to do so using suitably deformed Hermitian fiber
metrics. For gauge fields one can then rely on
\cite{jurco.et.al.:2001a, jurco.schupp.wess:2000a}, formulated in a
suitable geometric fashion.

In a third step, one wants to construct quantum field theories
corresponding to the classical field theories on a locally
noncommutative space-time. We shall discuss this in more detail in the
next section.

Note however, that for a consistent dynamical treatment we have to go
at least one step further: the Poisson structure $\theta$ itself (and
hence the corresponding star product) should be considered as a
dynamical quantity instead of a fixed background field. This is of
course desirable in any model of noncommutative space-time and
therefore poses a general open problem in such approaches to Planck
scale physics. Here, the principle difficulty is to impose field
equations on $\theta$ which are compatible with the Jacobi identity
$\Schouten{\theta, \theta} = 0$ and still allow some interesting
coupling to other fields.

%
% QFT
%

\subsection{Quantum Field Theories}
\label{subsec:QFT}

We continue with some heuristic remarks on possible quantum field
theories on a locally noncommutative space-time.

The main goal of our construction is the avoidance of a violation of
locality at large distances in field theory. To get a first glimpse on
how powerful our ideas might turn out to be, consider the flat
Minkowski space $M$ with a star product on $M \times M$ meeting our
requirements of compact support in relative coordinate directions.
Let us assume that the free field $\phi(x)$ is the ordinary one, i.e.
$\phi(x) = \int \hat \phi(k)\ \E^{ikx} \D k$ with the operator valued
distribution $\hat \phi(k)$ acting on Fock space as annihilation and
creation operators, respectively.

Now consider two fields at different points $x$ and $y$ in space time
(to be precise, we have to evaluate $\phi$ in suitable test functions
supported around $x$ and $y$). In order to fit them into our
framework, we consider the operator-valued distributions
\begin{equation}
    \label{eq:fgDef}
    f(x,y) = \phi(x)  \qquad \mbox{ and } \qquad g(x,y)=\phi(y)
\end{equation}
Using the globally defined exponential map, we now define for $p =
\frac{x+y}{2}$ and $v = \frac{x-y}{2}$
\begin{equation}
    \label{eq:FGDef}
    F(p,v)=  \phi(p+v)
    \qquad \mbox{ and } \qquad 
    G(p,v) = \phi(p-v)
\end{equation}
and consider their star product (ignoring all problems that arise in
taking such products of distributions). We then find that for $v
\notin \supp \theta_p$, the star product becomes the ordinary one and
in particular, the commutator of fields is unchanged,
\begin{equation}
    \label{eq:FieldCommutator}
    [F,G]_{\star_p} (p,v) 
    = [\phi(p+v), \phi(p-v)]_{\starp} = [\phi(p+v),\phi(p-v)] 
    \qquad \mbox{ for } \qquad v \notin \supp \theta_p
\end{equation}
Clearly, for $v$ large enough we thus recover micro-locality (in
particular, vanishing of the commutator for space-like separated
points). The classic no-go theorems on nonlocal fields
\cite{borchers.pohlmeyer:1968a, pohlmeyer:1968a} are circumvented as
the product of the fields itself is changed.

It remains an open and difficult problem to define a sensible
interaction term.  Possibly, one should employ the deformed
$\delta$-distributions (much in the spirit of
\cite{bahns.doplicher.fredenhagen.piacitelli:2003a}) or try to
generalize the approach to more than two points.

Our hope is however, that once this has been achieved, the field
theory's properties regarding renormalization should be considerably
improved compared to both the ordinary one as well as the one based on
models with constant $\theta$ (as the infrared regime is clearly
separated from the ultraviolet one). In the long run, it would have to
be investigated whether a quickly decreasing (instead of a compactly
supported) noncommutative structure would suffice for the purposes of
renormalization such that problems with the Lorentz structure could be
avoided -- though the construction possibly only admits this for flat
space (where the exponential map defines a global diffeomorphism).

%
% Further extensions of the model
%

\subsection{Further extensions of the model}
\label{subsec:FurtherExtensions}

Let us finally mention some possible extensions of the locally
noncommutative space-times as presented above. Once having realized
that not $M$ but $M \times M$ is relevant when discussing small
distance behavior one can of course go one step further: In principle
one can also discuss noncommutativity which only becomes present when
\emph{three} points come close together. This would give a
noncommutativity on $M \times M \times M$ located again close to the
diagonal. Analogously, one can consider Poisson structures
$\theta^{(k)}$ on $M^k = M \times \cdots \times M$ for arbitrary $k
\in \mathbb{N}$ with support close around the total diagonal. It would
be clearly a very interesting investigation how one can combine all
these $\theta^{(k)}$ and formulate compatibilities between them for
different $k$. In particular, it would be interesting to find a
reasonable replacement for the verticality requirement.

The last extension we want to mention is the passage from formal star
products to convergent deformations. This is known to be a serious
problem in deformation quantization and not much can be said on a
general level. However, for certain Poisson structures there exist
convergent star products quantizing them, usually by means of suitable
integral formulas, see e.g.~\cite{bieliavsky:2002a, rieffel:1993a,
  landsman:1998a}.  Eventually, the result will be a $C^*$-algebraic
approach like in \cite{doplicher.fredenhagen.roberts:1995a} which will
be necessary for all questions concerning quantum field theories in
the locally noncommutative space-time. For vertical Poisson structures
$\theta$ arising from actions of some $\mathbb{R}^d$ one can rely on
Rieffel's general construction \cite{rieffel:1993a} to obtain a
$C^*$-algebraic deformation. This will be investigated in a future
project \cite{heller.neumaier.waldmann:2006a}.

%
% Appendix: vertical formalities
%

\appendix
\section{Vertical formality and vertical star products}
\label{sec:VerticalFormality}

In this appendix we collect some results on star products on vector
bundles which seem to be new but follow essentially in a
straightforward manner from Kontsevich's formality theorem for
$\mathbb{R}^d$. Thus we only indicate the proofs and outline the
ideas.

%
% Vertical multivector fields
%

\subsection{Vertical multivector fields on a vector bundle}
\label{subsec:VerticalMultivectors}

Let us first recall some standard results in order to fix our
notation. We consider a real vector bundle $\pi: E \longrightarrow M$
with with fiber dimension $d$.

Recall that a tangent vector $X \in T_vE$ at $v \in E$ is called
\emph{vertical} if $T_v\pi X = 0$. The subbundle of vertical tangent
vectors is denoted by $\Ver(E) \subseteq TE$. Moreover, we obtain
vertical contravariant tensors $\bigotimes^k \Ver(E) \subseteq
\bigotimes^k TE$ as well, in particular the vertical symmetric and
antisymmetric contravariant tensor bundles $\Sym^k \Ver(E) \subseteq
\Sym^k TE$ and $\Anti^k\Ver(E) \subseteq \Anti^k TE$, respectively.
The corresponding sections are the vertical contravariant tensor
fields $\Gamma^\infty(\bigotimes^k \Ver(E)) \subseteq
\Gamma^\infty(\bigotimes^kTE)$ where we are most interested in the
vertical multivector fields $\Xver^\bullet(E) =
\Gamma^\infty(\Anti^\bullet\Ver(E))$.

In the sequel we make use of local expressions. Thus fix a locally
defined basis of sections $e_1, \ldots, e_d \in \Gamma^\infty(E|_U)$
where $U \subseteq M$ is a suitable open subset and $E|_U =
\pi^{-1}(U)$. The corresponding dual basis is denoted by $e^1, \ldots,
e^d \in \Gamma^\infty(E^*|_U)$. The choice of such a basis induces
\emph{linear coordinates} $s^1, \ldots, s^d \in C^\infty(E|_U)$ along
the fibers by setting $s^\alpha (v) = \SP{e^\alpha(\pi(v)), v}$ as
usual, locally trivializing $E|_U \cong U \times \mathbb{R}^d$. If
$(x^1, \ldots, x^n)$ are local coordinates on $U \subseteq M$ then
$(x^1 \circ \pi, \ldots, x^n \circ \pi, s^1, \ldots, s^d)$ are local
coordinates of $E$ defined on $E|_U$. The local tangent vector fields
$\frac{\partial}{\partial s^\alpha}$ are vertical and provide a basis
of sections for $\Ver(E)|_U$.

We can \emph{lift} sections $s \in \Gamma^\infty(E)$ to vertical
vector fields $s^\ver \in \Xver^1(E)$ by setting
\begin{equation}
    \label{eq:VerticalLift}
    s^\ver (v) = 
    \frac{\D}{\D t} \Big|_{t=0} \left(v + t s(\pi(v))\right)
\end{equation}
for $v \in E$, whence clearly $\frac{\partial}{\partial s^\alpha} =
e_\alpha^\ver$. Hence the tangent vector fields
$\frac{\partial}{\partial s^\alpha}$ do not depend on the choice of
the local coordinates $(x^1, \ldots, x^d)$ but only on the frame $e_1,
\ldots, e_d$.

We can extend the vertical lift to arbitrary tensor fields $X \in
\Gamma^\infty(\bigotimes^kE)$ in the usual way, compatible with the
tensor product, where a $0$-tensor field $u \in
\Gamma^\infty(\bigotimes^0E) = C^\infty(M)$ is lifted via $u^\ver =
\pi^*u$.  Locally, any vertical tensor field $X \in
\Gamma^\infty(\bigotimes^k \Ver(E))$ can be written as
\begin{equation}
    \label{eq:XVerlocal}
    X\big|_{\pi^{-1}(U)} =
    X^{\alpha_1 \cdots \alpha_k} 
    e_{\alpha_1}^\ver \otimes \cdots \otimes e_{\alpha_k}^\ver
\end{equation}
with $X^{\alpha_1 \cdots \alpha_k} \in C^\infty(\pi^{-1}(U))$. Then
$X$ is a vertical lift iff the functions $X^{\alpha_1 \cdots
  \alpha_k}$ are pull-backs of functions in $C^\infty(U)$.  Denote by
$\xi \in \Gamma^\infty(\Ver(E))$ the \emph{Euler vector field},
defined via its flow $(t, v) \mapsto \E^t v$.  Locally,
$\xi|_{\pi^{-1}(U)} = s^\alpha \frac{\partial}{\partial s^\alpha}$ and
hence $X \in \Gamma^\infty(\bigotimes^k \Ver(E))$ is a vertical lift
iff $\Lie_\xi X = -kX$.

More generally, we say that $X \in \Gamma^\infty(\bigotimes^k
\Ver(E))$ is \emph{polynomial along the fibers} of degree $\ell$ if
$\Lie_\xi X = (\ell-k) X$. Clearly, this is equivalent to the local
statement that all coefficient functions $X^{\alpha_1 \cdots
  \alpha_k}$ are polynomials in the fiber variables $s^1, \ldots, s^d$
of degree $\ell$. The vertical tensor fields $X \in
\Gamma^\infty(\bigotimes^k \Ver(E))$ polynomial along the fibers of
degree $\ell$ are in canonical bijection to tensor fields
$\widetilde{X} \in \Gamma^\infty(\Sym^\ell E^* \otimes \bigotimes^k
E)$ via the relation
\begin{equation}
    \label{eq:XverpolyTildeX}
    X(v) = \left(\widetilde{X}\big|_{\pi(v)}(v, \ldots, v)\right)^\ver,
\end{equation}
where we first insert the point $v \in E_{\pi(v)}$ in the $\Sym^\ell
E^*$-part $\ell$-times and then lift the $\bigotimes^k E$-part
vertically. In particular, vertical lifts are those vertical tensor
fields which are \emph{constant along the fibers}.

A vertical tensor field $X \in \Gamma^\infty(\bigotimes^k\Ver(E))$ can
be \emph{restricted} to a fiber $E_p \subseteq E$ for $p \in M$ and
yields a tensor field in $\Gamma^\infty(\bigotimes^k TE_p)$. This
follows from the fact that canonically $\ker T_v\pi \cong T_v(E_p)$
for $v \in E_p$. Let $\iota_p: E_p \hookrightarrow E$ denote the
inclusion map then we denote the restriction by $\iota_p^* X \in
\Gamma^\infty(\bigotimes^k TE_p)$. In particular, if $X \in
\Gamma^\infty(\bigotimes^k \Ver(E))$ is polynomial along the fibers of
degree $\ell$ then $\iota_p^* X$ is a tensor field on the \emph{vector
  space} $E_p$ which is polynomial of degree $\ell$ in the usual
sense.

Finally, we focus on vertical multivector fields. The following is
folklore and consists in a straightforward verification:
\begin{proposition}
    \label{proposition:VerticalSubGerstenhaber}
    Let $\pi: E \longrightarrow M$ be a vector bundle.
    \begin{enumerate}
    \item The vertical multivector fields $\Xver^\bullet(E)$ are a
        Gerstenhaber subalgebra of $\mathfrak{X}^\bullet(E)$.
    \item The restriction map $\iota_p^*: \Xver^\bullet(E)
        \longrightarrow \mathfrak{X}^\bullet(E_p)$ is a surjective
        homomorphism of Gerstenhaber algebras, i.e. for all $X, Y \in
        \Xver^\bullet(E)$ we have
        \begin{equation}
            \label{eq:RestrictGerstenhaber}
            \iota_p^*(X \wedge Y) = \iota_p^* X \wedge \iota_p^* Y
            \quad
            \textrm{and}
            \quad
            \iota_p^*(\Schouten{X, Y}) 
            = \Schouten{\iota_p^*X, \iota_p^* Y}.
        \end{equation}
    \item The vertical multivector fields which are polynomial along
        the fibers are a Gerstenhaber subalgebra of $\Xver^\bullet(E)$
        isomorphic to the Gerstenhaber algebra
        $\bigoplus_{\ell=0}^\infty \Gamma^\infty(\Sym^\ell E^* \otimes
        \Anti^\bullet E)$, equipped with its canonical fiberwise
        Gerstenhaber algebra structure.
    \end{enumerate}
\end{proposition}
Note that for $X, Y \in \Gamma^\infty(\Anti^\bullet TE)$ we have
\begin{equation}
    \label{eq:VerticalLiftsCommute}
    \Schouten{X^\ver, Y^\ver} = 0.
\end{equation}

%
% The vertical Hochschild-Kostant-Rosenberg theorem
%

\subsection{The vertical Hochschild-Kostant-Rosenberg theorem}
\label{subsec:VerticalHKR}

Recall that $k$-vector fields $X \in \mathfrak{X}^k(E)$ can be viewed
as totally antisymmetric first order $k$-differential operators by use
of the \emph{Hochschild-Kostant-Rosenberg map} (for short: HKR map)
\begin{equation}
    \label{eq:HKRDef}
    \left(\mathrm{U}^{(1)}(X)\right)(f_1, \ldots, f_k)
    =
    \frac{1}{k!} \SP{X, \D f_1 \otimes \cdots \otimes \D f_k},
\end{equation}
where $f_1, \ldots, f_k \in C^\infty(E)$. We denote by
$\HCdiff^k(C^\infty(E))$ the differential Hochschild $k$-cochains with
values in $C^\infty(E)$, i.e. the $k$-differential operators
\begin{equation}
    \label{eq:kcochains}
    \phi: \underbrace{C^\infty(E) \times \cdots \times
      C^\infty(E)}_{k\textrm{-times}}
    \longrightarrow C^\infty(E).
\end{equation}
Then $\phi \in \HCdiff^k(C^\infty(E))$ is called \emph{vertical} if
\begin{equation}
    \label{eq:VerticalMultiDiffOp}
    \phi(f_1, \ldots, \pi^*u f_i, \ldots, f_k)
    =
    \pi^*u \: \phi(f_1, \ldots, f_k)
\end{equation}
for all $f_1, \ldots, f_k \in C^\infty(E)$, $u \in C^\infty(M)$ and $i
= 1, \ldots, k$. We denote the vertical $k$-differential operators by
$\HCdiffver^k(C^\infty(E))$. Clearly, $\mathrm{U}^{(1)} (X) \in
\HCdiffver^k(C^\infty(E))$ for a vertical $k$-vector field $X$. The
restriction of $\mathrm{U}^{(1)}$ to vertical multivector fields is
denoted by
\begin{equation}
    \label{eq:HKRver}
    \Uver^{(1)}: \Xver^\bullet(E) 
    \longrightarrow \HCdiffver^\bullet(C^\infty(E)).
\end{equation}

Let $R = (r_1, \ldots, r_k) \in \mathbb{N}^k$ be the multi-order of $\phi
\in \HCdiffver^k(C^\infty(E))$. Then locally
\begin{equation}
    \label{eq:DverLocal}
    \phi(f_1, \ldots, f_k) \big|_{\pi^{-1}(U)}
    =
    \sum_{L=0}^R \frac{1}{L!}
    \phi_L^{\alpha_1^1 \cdots \alpha_{\ell_1}^1 
      \cdots \alpha_1^k \cdots \alpha_{\ell_k}^k}
    \frac{\partial^{\ell_1} f_1}
    {\partial s^{\alpha_1^1} \cdots \partial s^{\alpha_{\ell_1}^1}}
    \cdots
    \frac{\partial^{\ell_k} f_k}
    {\partial s^{\alpha_1^k} \cdots \partial s^{\alpha_{\ell_k}^k}}
\end{equation}
with unique functions $\phi_L^{\alpha_1^1 \cdots \alpha_{\ell_1}^1
  \cdots \alpha_1^k \cdots \alpha_{\ell_k}^k} \in
C^\infty(\pi^{-1}(U))$, totally symmetric in $(\alpha_1^i, \ldots,
\alpha_{\ell_i}^i)$ for all $i = 1, \ldots, k$. Conversely, if $\phi
\in \HCdiff^k(C^\infty(E))$ is locally of this form, then $\phi$ is
vertical. In this case it is easy to see that for all $L \le R$
\begin{equation}
    \label{eq:SymbolLphi}
    \sigma_L(\phi) \big|_{\pi^{-1}(U)}
    =
    \frac{1}{L!}
    \phi_L^{\alpha_1^1 \cdots \alpha_{\ell_1}^1 
      \cdots \alpha_1^k \cdots \alpha_{\ell_k}^k}
    \frac{\partial}{\partial s^{\alpha_1^1}} 
    \vee \cdots \vee
    \frac{\partial}{\partial s^{\alpha_{\ell_1}^1}}
    \otimes \cdots \otimes
     \frac{\partial}{\partial s^{\alpha_1^k}} 
    \vee \cdots \vee
    \frac{\partial}{\partial s^{\alpha_{\ell_k}^k}}
\end{equation}
defines a global tensor field $\sigma_L(\phi) \in
\Gamma^\infty(\Sym^{\ell_1} \Ver(E) \otimes \cdots \otimes
\Sym^{\ell_k} \Ver(E))$, the \emph{$L$-symbol} of $\phi$. Note that in
general only the leading symbol, i.e. for $L = R$, has an invariant
geometric interpretation as a tensor field. Conversely, out of a given
tensor field $A \in \Gamma^\infty(\Sym^{\ell_1} \Ver(E) \otimes \cdots
\otimes \Sym^{\ell_k} \Ver(E))$ one can build a unique vertical
$k$-differential operator $Q_L(A)$ of multi-order $L$ in a canonical
way, such that $\sigma_L(Q_L(A)) = A$ and $\sigma_{L'}(Q_L(A)) = 0$
for all $L' \neq L$.

Vertical multidifferential operators can again be restricted to fibers
and yield multidifferential operators on $E_p$ for each $p \in M$. We
denote the restriction again by
\begin{equation}
    \label{eq:iotap}
    \iota_p^*: 
    \HCdiffver^\bullet(C^\infty(E))
    \longrightarrow
    \HCdiff^\bullet(C^\infty(E_p)).
\end{equation}
Vertical multidifferential operators behave well under
multicomposition: if $\phi_i \in \HCdiffver^{\ell_i}(C^\infty(E))$ for
$i = 1, \ldots, k$ and $\phi \in \HCdiffver^k(C^\infty(E))$ are given
then $\phi \circ (\phi_1 \otimes \cdots \otimes \phi_k) \in
\HCdiffver^{\ell_1 + \cdots + \ell_k}(C^\infty(E))$. This is obvious
from the definition \eqref{eq:VerticalMultiDiffOp}. Moreover, in this
case we have
\begin{equation}
    \label{eq:RestrictComposition}
    \iota^*_p \left(
        \phi \circ (\phi_1 \otimes \cdot \otimes \phi_k)
    \right)
    =
    \left(\iota^*_p \phi\right) \circ 
    \left(
        \iota^*_p \phi_1 \otimes \cdots \otimes \iota^*_p \phi_k
    \right).
\end{equation}
From this we immediately have the following statement:
\begin{proposition}
    \label{proposition:VerticalHC}
    Let $\pi: E \longrightarrow M$ be a vector bundle and $p \in M$.
    \begin{enumerate}
    \item The vertical differential Hochschild cochains
        $\HCdiffver^\bullet(C^\infty(E))$ are closed under the
        cup-product $\cup$, the Hochschild differential $\delta$, and
        the Gerstenhaber bracket $[\cdot, \cdot]$.
    \item The restriction map
        \begin{equation}
            \label{eq:iotapHCMorph}
            \iota^*_p: \HCdiffver^\bullet(C^\infty(E))
            \longrightarrow \HCdiff^\bullet(C^\infty(E_p))
        \end{equation}
        is compatible with the cup-products, the Hochschild
        differentials and the Gerstenhaber brackets, respectively.
    \end{enumerate}
\end{proposition}
\begin{proof}
    See e.g.~\cite{gerstenhaber:1963a} for the definition of $\cup$,
    $\delta$, and $[\cdot, \cdot]$. The statement follows immediately
    from the compatibility with the multicomposition and
    \eqref{eq:RestrictComposition}.
\end{proof}

In particular, $\HCdiffver^\bullet(C^\infty(E))$ is a subcomplex of
the usual differential Hochschild complex of $C^\infty(E)$. Its
cohomology, the \emph{vertical Hochschild cohomology} of
$C^\infty(E)$, will be denoted by $\HHdiffver^\bullet(C^\infty(E))$.
It is well-known that the cup-product and the Gerstenhaber bracket
pass to the differential Hochschild cohomology
$\HHdiff^\bullet(C^\infty(E))$ which then becomes a Gerstenhaber
algebra. To show the appropriate algebraic identities between $\cup$
and $[\cdot, \cdot]$ one has to construct explicit coboundaries, see
\cite{gerstenhaber:1963a}. From these explicit formulas it can easily
be seen that the relevant coboundaries can be chosen vertical if all
involved cocycles are vertical. Hence one has the following result:
\begin{proposition}
    \label{proposition:HHdiffverRestriction}
    Let $\pi: E \longrightarrow M$ be a vector bundle and $p \in M$.
    \begin{enumerate}
    \item The vertical differential Hochschild cohomology of
        $C^\infty(E)$ becomes a Gerstenhaber algebra with respect to
        the cup-product and the Gerstenhaber bracket. The canonical
        map
        \begin{equation}
            \label{eq:IncludeHHverHH}
            \HHdiffver^\bullet(C^\infty(E)) 
            \longrightarrow 
            \HHdiff^\bullet(C^\infty(E))
        \end{equation}
        is a map of Gerstenhaber algebras.
    \item The restriction map $\iota_p^*$ induces a map of
        Gerstenhaber algebras
        \begin{equation}
            \label{eq:RestrictHHver}
            \iota_p^*: \HHdiffver^\bullet(C^\infty(E))
            \longrightarrow
            \HHdiff^\bullet(C^\infty(E_p)).
        \end{equation}
    \end{enumerate}
\end{proposition}

Let us now compute $\HHdiffver^\bullet(C^\infty(E))$. We start with
the trivial vector bundle $E = U \times \mathbb{R}^d$. In this case we
consider $C^\infty(U \times \mathbb{R}^d)$ as a symmetric
$C^\infty(\mathbb{R}^d)$-bimodule in the canonical way. If $\phi \in
\HCdiff^k(C^\infty(\mathbb{R}^d), C^\infty(U \times
\mathbb{R}^d))$ is a $k$-differential operator on
$C^\infty(\mathbb{R}^d)$ with values in $C^\infty(U \times
\mathbb{R}^d)$ we can view this as a vertical $k$-differential
operator $\widetilde{\phi} \in \HCdiffver^k(C^\infty(U \times
\mathbb{R}^d))$ by setting
\begin{equation}
    \label{eq:Tildephi}
    \left(\widetilde{\phi} (f_1, \ldots, f_k) \right)(u, v)
    =
    \left(\phi(f_1(u, \cdot), \ldots, f_k(u, \cdot))\right)(u, v).
\end{equation}
Conversely, let $\Phi \in \HCdiffver^k(C^\infty(U \times
\mathbb{R}^d))$ then we can simply restrict $\Phi$ to the subalgebra
$C^\infty(\mathbb{R}^d)$ of $C^\infty(U \times \mathbb{R}^d)$ and
obtain $\widehat{\Phi} \in \HCdiff^k(C^\infty(\mathbb{R}^d),
C^\infty(U \times \mathbb{R}^d))$. The following is obvious:
\begin{lemma}
    \label{lemma:IsomorphicComplex}
    The extension $\phi \mapsto \widetilde{\phi}$ and the restriction
    $\Phi \mapsto \widehat{\Phi}$ are mutually inverse isomorphisms of
    Hochschild complexes
    \begin{equation}
        \label{eq:IsomorphicComplex}
        \HCdiff^k(C^\infty(\mathbb{R}^d), 
        C^\infty(U \times \mathbb{R}^d))
        \cong        
        \HCdiffver^k(C^\infty(U \times \mathbb{R}^d)).
    \end{equation}
\end{lemma}

Since $C^\infty(U \times \mathbb{R}^d)$ is a \emph{symmetric}
$C^\infty(\mathbb{R}^d)$-bimodule, we easily can compute its
differential Hochschild cohomology using the Koszul `resolution' as in
\cite{bordemann.et.al:2005a:pre}:
\begin{lemma}
    \label{lemma:Koszul}
    The HKR map induces an isomorphism
    \begin{equation}
        \label{eq:HKRIso}
        \Xver^\bullet(U \times \mathbb{R}^d)
        \cong
        \Anti^\bullet \mathbb{R}^d \otimes 
        C^\infty(U \times \mathbb{R}^d)
        \longrightarrow
        \HHdiff^\bullet (C^\infty(\mathbb{R}^d), 
        C^\infty(U \times \mathbb{R}^d))
        \cong
        \HHdiffver (C^\infty(U \times \mathbb{R}^d)).
    \end{equation}
\end{lemma}
In particular, if $\phi \in \HCdiffver^\bullet(C^\infty(U \times
\mathbb{R}^d))$ is a cocycle then
\begin{equation}
    \label{eq:LocalCocycleUeinsX}
    \phi = \mathrm{U}^{(1)}(X) + \delta \psi
\end{equation}
for some $\psi \in \HCdiffver^{\bullet -1}(C^\infty(U\times
\mathbb{R}^d))$ and a unique $X \in \Xver^\bullet(C^\infty(U \times
\mathbb{R}^d))$, given by the total antisymmetrization of $\phi$.

From this local statement the standard globalization argument as
e.g. in \cite{cahen.gutt.dewilde:1980a} using a suitable partition of
unity of $M$ and local trivializations of $E$ gives the following
statement:
\begin{lemma}
    \label{lemma:GlobalHKR}
    If $\phi \in \HCdiffver^\bullet(C^\infty(E))$ is a cocycle then
    there exists a unique vertical multivector field $X \in
    \Xver^\bullet(E)$ and a vertical $\psi \in
    \HCdiffver^{\bullet-1}(C^\infty(E))$ such that
    \begin{equation}
        \label{eq:GlobalCocyclePhiUeinsXPsi}
        \phi = \Uver^{(1)}(X) + \delta \psi,
    \end{equation}
    and $X$ is obtained from the total antisymmetrization of
    $\phi$. In particular,
    \begin{equation}
        \label{eq:UeinsVectorspaceIso}
        \Uver^{(1)}: \Xver^\bullet(E) 
        \longrightarrow
        \HHdiffver^\bullet(C^\infty(E))
    \end{equation}
    is an isomorphism.
\end{lemma}

Finally, $\Uver^{(1)}$ as in \eqref{eq:UeinsVectorspaceIso} is not
only an isomorphism of vector spaces but compatible with the
Gerstenhaber algebra structures: First it is clear that on the level
of cochains $\Uver^{(1)}$ maps $\wedge$-products to the total
antisymmetrization of the corresponding $\cup$-products. Passing to
cohomology, $\cup$ becomes supercommutative whence $\Uver^{(1)}$ maps
$\wedge$-products to $\cup$-products in cohomology. Moreover,
$\Uver^{(1)}$ is easily verified to map Schouten brackets of functions
and vector fields to the corresponding Gerstenhaber brackets already
on the level of cochains. Since functions and vector fields generate
$\Xver^\bullet(E)$ by $\wedge$-products and since $\Schouten{\cdot,
  \cdot}$ as well as $[\cdot, \cdot]$ satisfy the same Leibniz rule
(the latter only in cohomology) and since $\Uver^{(1)}$ is an
isomorphism of associative supercommutative algebras, it follows that
$\Uver^{(1)}$ also maps Schouten brackets to Gerstenhaber brackets in
cohomology. Note that on the level of cochains this is not true for
higher multivector fields. We summarize the result of this section:
\begin{theorem}[Vertical Hochschild-Kostant-Rosenberg Theorem]
    \label{theorem:HKRver}
    Let $\pi: E \longrightarrow M$ be a vector bundle and $p \in M$.
    \begin{enumerate}
    \item The vertical HKR map gives an isomorphism of Gerstenhaber algebras
        \begin{equation}
            \label{eq:verHKRIsoGerstenhaber}
            \Uver^{(1)}: \Xver^\bullet(E) 
            \longrightarrow 
            \HHdiffver^\bullet(C^\infty(E)).
        \end{equation}
    \item Let $\mathrm{U}_p^{(1)}: \mathfrak{X}^\bullet(E_p)
        \longrightarrow \HHdiff^\bullet(C^\infty(E_p))$ be the usual
        HKR map on the vector space $E_p$. Then
        \begin{equation}
            \label{eq:HKRCommute}
            \bfig
            \square<1000,400>[\Xver^\bullet(E)%
            `\HCdiffver^\bullet(C^\infty(E))%
            `\mathfrak{X}^\bullet(E_p)%
            `\HCdiff(C^\infty(E_p))%
            ;\Uver^{(1)}%
            `\iota_p^*%
            `\iota_p^*%
            `\mathrm{U}^{(1)}_p]%
            \efig
        \end{equation}
        commutes and all maps are homomorphisms of Gerstenhaber
        algebras.
    \end{enumerate}
\end{theorem}

%
% The vertical formality
%

\subsection{The vertical formality}
\label{subsec:VerticalFormality}

We come now to the main theorem of this appendix, for which we shall
recall some basic notions of formal deformation theory, see e.g.
\cite{cattaneo:2005a} or \cite[Appendix]{bordemann.et.al:2005a:pre},
and the language of coalgebras, see
e.g.~\cite{markl.shnider.stasheff:2002a, sweedler:1969a}. Let
$\mathfrak{g} = \bigoplus_{k \in \mathbb{Z}} \mathfrak{g}_k$ be a
differential graded Lie algebra with Lie bracket $[\cdot, \cdot]$ and
differential $\delta$ of degree $+1$. This structure can alternatively
be described as follows. We denote the same vector space with shifted
degree by $+1$ by $\mathfrak{g}[1]$ and consider the graded symmetric
algebra $\Sym(\mathfrak{g}[1])$. With the graded symmetric tensor
product $\vee$ and the graded cocommutative shuffle coproduct
$\Delta_{\mathrm{sh}}$ one obtains a bialgebra $\Sym(\mathfrak{g}[1])$
with unit $\Unit$ and counit $\epsilon$ being just the projection on
the tensor degree $0$. As coalgebra, $\Sym(\mathfrak{g}[1])$ is cofree
within the category of augmented graded cocommutative counital
coalgebras $\mathcal{CC}_{AN}$ with nilpotent augmentation ideal $\ker
\epsilon$, where a coalgebra is called \emph{augmented} if there is
exactly one group-like element denoted by $\Unit$.  Note that
$\Sym(\mathfrak{g}[1])$ is \emph{not} cofree within the category of
all graded cocommutative coalgebras, see
\cite{markl.shnider.stasheff:2002a, sweedler:1969a,
  bordemann.et.al:2005a:pre} for further details.  The differential
$\delta$ and the bracket $[\cdot, \cdot]$ can be combined to a single
map $\D = \delta + [\cdot, \cdot]: \Sym(\mathfrak{g}[1])
\longrightarrow \mathfrak{g}[1]$ of degree $+1$. Since
$\Sym(\mathfrak{g}[1])$ is in $\mathcal{CC}_{AN}$ and thanks to the
cofreeness, this map extends uniquely to a coderivation
$\underline{\D}: \Sym(\mathfrak{g}[1]) \longrightarrow
\Sym(\mathfrak{g}[1])$ such that $\pr_{\mathfrak{g}[1]} \circ
\underline{\D} = \D$. Then $\delta^2 = 0$, the compatibility between
$\delta$ and $[\cdot, \cdot]$ and the Jacobi identity for $[\cdot,
\cdot]$ are all encoded in $\underline{\D}^2 = 0$.

Generalizing this gives the definition of an \emph{$L_\infty$-algebra}
(or \emph{Lie algebra up to homotopy}): a graded vector space
$\mathfrak{g}$ is called $L_\infty$-algebra if $\Sym(\mathfrak{g}[1])$
is equipped with a coderivation $\underline{\D}$ of degree $+1$ and
$\underline{\D}^2 = 0$. If $\mathfrak{g}$ and $\mathfrak{h}$ are
$L_\infty$-algebras then an $L_\infty$-morphism is a coalgebra
morphism $\underline{\mathrm{U}}: \Sym(\mathfrak{g}[1])
\longrightarrow \Sym(\mathfrak{h}[1])$ such that
$\underline{\mathrm{U}} \circ \underline{\D}_{\mathfrak{g}} =
\underline{\D}_{\mathfrak{h}} \circ \underline{\mathrm{U}}$. Clearly,
in the above example of a differential graded Lie algebra, any
morphism of differential graded Lie algebras induces an
$L_\infty$-morphism. \emph{However}, in general, there are more
general $L_\infty$-morphisms than these: this is the key idea of
Kontsevich's formality theorem.

If $\mathfrak{g}$ has an $L_\infty$-structure then the coderivation
$\underline{\D}$ is uniquely determined by $\D = \pr_{\mathfrak{g}[1]}
\circ \underline{\D}$. Similarly, an $L_\infty$-morphism
$\underline{\mathrm{U}}$ is determined by $\mathrm{U} =
\pr_{\mathfrak{h[1]}} \circ \underline{U}$. Finally, each $d$ and
$\mathrm{U}$ are determined by their \emph{Taylor coefficients}
\begin{equation}
    \label{eq:TaylorCoefficients}
    \D = \sum_{n=1}^\infty \D^{(n)}
    \quad
    \textrm{and}
    \quad
    \mathrm{U} = \sum_{n=0}^\infty \mathrm{U}^{(n)},
\end{equation}
where $\D^{(n)}: \Sym^n(\mathfrak{g}[1]) \longrightarrow
\mathfrak{g}[1]$ and $\mathrm{U}^{(n)}: \Sym^n(\mathfrak{g}[1])
\longrightarrow \mathfrak{h}[1]$, respectively. Necessarily,
$\mathrm{U}^{(0)}$ maps $\Unit_{\Sym(\mathfrak{g}[1])}$ to
$\Unit_{\Sym(\mathfrak{h}[1])}$.

In general, $\delta = \D^{(1)}: \mathfrak{g}[1] \longrightarrow
\mathfrak{g}[1]$ satisfies $\delta^2 = 0$ and $[\cdot, \cdot] =
\D^{(2)}$ defines a Lie bracket `up to $\delta$-homotopy', i.e.
$[\cdot, \cdot]$ induces a graded Lie bracket on the
$\delta$-cohomology. This explains the name $L_\infty$-algebra.
Moreover, an $L_\infty$-morphism induces a morphism of graded Lie
algebras in cohomology. One calls $\underline{\mathrm{U}}$ an
\emph{$L_\infty$-quasiisomorphism} or \emph{formality} if the induced
map in cohomology is an isomorphism.

In \cite{kontsevich:2003a} Kontsevich constructed an
$L_\infty$-quasiisomorphism $\underline{\mathrm{U}}$ between
$\mathfrak{g} = \mathfrak{X}^\bullet(\mathbb{R}^d)[1]$ viewed as
$L_\infty$-algebra with $\D = \D^{(2)} = \Schouten{\cdot, \cdot}$, and
$\mathfrak{h} = \HCdiff^\bullet(C^\infty(\mathbb{R}^d))[1]$ viewed as
$L_\infty$-algebra with $\D = \D^{(1)} + \D^{(2)} = \delta + [\cdot,
\cdot]$. We recall the basic properties of this formality map:
\begin{theorem}[Kontsevich's Formality for $\mathbb{R}^d$]
    \label{theorem:KontsevichFormality}
    There exists an $L_\infty$-quasiisomorphism
    \begin{equation}
        \label{eq:KontsevichFormality}
        \underline{\mathrm{U}}_{\mathbb{R}^d}:
        \Sym(\mathfrak{X}^\bullet(\mathbb{R}^d)[2])
        \longrightarrow
        \Sym(\HCdiff^\bullet(C^\infty(\mathbb{R}^d))[2]),
    \end{equation}
    such that the Taylor coefficients
    $\mathrm{U}^{(n)}_{\mathbb{R}^d}$ have the following properties:
    \begin{enumerate}
    \item $\mathrm{U}^{(1)}_{\mathbb{R}^d}$ is the HKR map.
    \item $\mathrm{U}^{(n)}_{\mathbb{R}^d}$ is a real $n$-differential
        operator on its $n$ arguments in
        $\mathfrak{X}^\bullet(\mathbb{R}^d)$ with constant
        coefficients.
    \item $\mathrm{U}^{(n)}_{\mathbb{R}^d}$ is $\mathrm{GL}(d,
        \mathbb{R})$-invariant in the sense that for $X_1, \ldots, X_n
        \in \mathfrak{X}^\bullet(\mathbb{R}^d)$, $f_1, \ldots, f_m \in
        C^\infty(\mathbb{R}^d)$, and $A \in \mathrm{GL}(d,
        \mathbb{R})$ one has
        \begin{equation}
            \label{eq:GLdCovariance}
            A^* 
            \left(
                \big(
                \mathrm{U}^{(n)}_{\mathbb{R}^d} (X_1, \ldots, X_n)
                \big)
                (f_1, \ldots, f_m)
            \right)
            =
            \big(
            \mathrm{U}^{(n)}_{\mathbb{R}^d}(A^* X_1, \ldots, A^* X_n)
            \big)
            (A^* f_1, \ldots, A^* f_m),
        \end{equation}
        where $A$ acts by $x \mapsto Ax$ as usual.
    \end{enumerate}
\end{theorem}

The important point for us is the
$\mathrm{GL}(d,\mathbb{R})$-invariance. Using a local frame $e_1,
\ldots, e_d \in \Gamma^\infty(E|_U)$, on a local trivialization $E|_U
\cong U \times \mathbb{R}^d$ we can define
\begin{equation}
    \label{eq:UverDef}
    \Uver^{(n)}(X_1, \ldots, X_n) \big|_{E_p}
    =
    \mathrm{U}^{(n)}_{\mathbb{R}^d} 
    \left(X_1(u, \cdot), \ldots, X_n(u, \cdot)\right),
\end{equation}
for $X_1, \ldots, X_n \in \Xver^\bullet(E)$ and $p \in U$. We only use
the linear coordinates $(s^1, \ldots, s^d)$ on $E|_U$ and apply
$\mathrm{U}^{(n)}_{\mathbb{R}^d}$ with respect to those, treating the
$U$-directions as parameters not affected by
$\mathrm{U}^{(n)}_{\mathbb{R}^d}$. Then from \eqref{eq:GLdCovariance}
it follows immediately, that for vertical $X_1, \ldots, X_n$ the
operator $\Uver^{(n)}(X_1, \ldots, X_n)$ is actually defined globally
and independent on the choice of the trivialization.  From this we
obtain immediately the main result of this appendix, the vertical
formality theorem:
\begin{theorem}[Vertical formality theorem]
    \label{theorem:VerticalFormality}
    Let $\pi: E \longrightarrow M$ be a vector bundle. Then there
    exists a unique $L_\infty$-quasiisomorphism
    \begin{equation}
        \label{eq:Uver}
        \uUver:
        \Sym(\Xver^\bullet(E)[2]) 
        \longrightarrow 
        \Sym(\HCdiffver^\bullet(C^\infty(E))[2])
    \end{equation}
    which has the following properties:
    \begin{enumerate}
    \item The Taylor coefficients $\Uver^{(n)}$ are real vertical
        $n$-differential operators on their $n$ arguments with
        constant coefficients.
    \item $\Uver^{(1)}$ is the HKR map.
    \item $\uUver$ restricts to an $L_\infty$-quasiisomorphism
        $\underline{\mathrm{U}}_p$
        \begin{equation}
            \label{eq:UverRestrict}
            \bfig
            \square<1200,400>[\Sym(\Xver^\bullet(E)\lbrack2\rbrack)%
            `\Sym(\HCdiffver^\bullet(C^\infty(E))\lbrack2\rbrack)%
            `\Sym(\mathfrak{X}^\bullet(E_p)\lbrack2\rbrack)%
            `\Sym(\HCdiff^\bullet(C^\infty(E_p))\lbrack2\rbrack)%
            ;\uUver%
            `\underline{\iota}^*_p%
            `\underline{\iota}^*_p%
            `\underline{\mathrm{U}}_p]%
            \efig
        \end{equation}
        such that $\underline{\mathrm{U}}_p$ coincides with
        Kontsevich's formality on the vector space $E_p$ for all $p
        \in M$. Here $\underline{\iota}^*_p$ is the canonical
        extension of the restriction $\iota_p^*$ to a coalgebra
        morphism.
    \end{enumerate}
\end{theorem}
\begin{proof}
    If $\uUver$ is a vertical $L_\infty$-morphism like in
    \eqref{eq:Uver} which satisfies the first part, then it clearly
    restricts to a $L_\infty$-morphism $\underline{\mathrm{U}}_p$ for
    all $p \in M$ such that \eqref{eq:UverRestrict} commutes.
    Moreover, such a $\uUver$ is completely determined by its
    restrictions $\underline{\mathrm{U}}_p$ which proves uniqueness.
    For existence, we see that $\uUver$ constructed above is an
    $L_\infty$-morphism since this can be checked locally whence we
    can rely on Theorem~\ref{theorem:KontsevichFormality}. Moreover,
    the first and third part are clearly satisfied by construction and
    the second part follows from
    Theorem~\ref{theorem:KontsevichFormality} as well as
    Theorem~\ref{theorem:HKRver}. Since by
    Theorem~\ref{theorem:HKRver} the vertical HKR map induces an
    isomorphism in cohomology, $\uUver$ is a quasiisomorphism since
    this is always decided by the first Taylor coefficient.
\end{proof}

%
% Vertical Poisson structures and vertical star products
%

\subsection{Vertical Poisson structures and vertical star products}
\label{subsec:VerticalPoissonStar}

A \emph{vertical Poisson structure} $\theta$ is a vertical bivector
field $\theta \in \Xver^2(E)$ with $\Schouten{\theta, \theta} =
0$. Analogously, one defines a \emph{formal vertical Poisson
  structure} $\theta = \sum_{r=0}^\infty \lambda^r \theta_r \in
\Xver^2(E)[[\lambda]]$. Two formal vertical Poisson structures
$\theta$, $\theta'$ are called \emph{vertically equivalent} if there
exists a formal vertical vector field $X \in \Xver^1(E)[[\lambda]]$
such that
\begin{equation}
    \label{eq:VerPoissonEquivalent}
    \theta' = \E^{\lambda\Lie_X} (\theta).
\end{equation}
One calls $\E^{\lambda\Lie_X}$ also a formal diffeomorphism.  Note that
in this case the zeroth order parts of $\theta$ and $\theta'$ coincide
$\theta_0 = \theta'_0$.

If $\theta$ is a vertical Poisson structure then $\theta_p = \iota^*_p
\theta \in \mathfrak{X}^2(E_p)$ is a Poisson structure on the vector
space $E_p$.  The map $\iota_p^*: C^\infty(E) \longrightarrow
C^\infty(E_p)$ then becomes a \emph{Poisson map}.  Clearly, vertically
equivalent $\theta, \theta' \in \Xver^2(E)[[\lambda]]$ restrict to
equivalent $\theta_p, \theta'_p \in \mathfrak{X}^2(E_p)[[\lambda]]$
via the restriction $X_p = \iota^*_p X$ of $X$.

A Poisson structure $\theta \in \mathfrak{X}^2(E)$ induces a Poisson
bracket $\{\cdot, \cdot\}_\theta$ on $C^\infty(E)$ as usual by $\{f,
g\}_\theta = \SP{\theta, \D f \otimes \D g}$. Then $\theta$ is vertical
iff $\{f, \pi^*u\}_\theta = 0$ for all $f \in C^\infty(E)$ and $u \in
C^\infty(M)$, i.e. $\pi^*C^\infty(M) \subseteq C^\infty(E)$ is part of
the Poisson center of $\{\cdot, \cdot\}_\theta$.

\begin{example}
    \label{example:VerticalConstant}
    Let $X \in \Gamma^\infty (\Anti^2 E)$ be an arbitrary
    section.  Then the vertical lift $\theta = X^\ver$ is a
    vertical Poisson structure by \eqref{eq:VerticalLiftsCommute}.
    Clearly, $\theta_p$ is a constant Poisson structure on $E_p$.
\end{example}
\begin{example}
    \label{example:BundleOfLieAlgebras}
    Let $E_p$ be equipped with a Lie algebra structure $[\cdot,
    \cdot]_p$ depending smoothly on $p \in M$. Then on $E^*$ we have a
    linear vertical Poisson structure in the usual way, which
    restricts to the canonical linear Poisson structure on each
    $E_p^*$ induced by the Lie bracket $[\cdot, \cdot]_p$.  More
    generally, one can consider quadratic and higher order vertical
    Poisson structures.
\end{example}
\begin{example}
    \label{example:VerticalAndCompactSupp}
    Consider first the local situation $E|_U \cong U \times
    \mathbb{R}^d$. Then we can choose $d$ commuting vector fields
    $X_1, \ldots, X_d$ on $\mathbb{R}^d$ whose supports are contained
    in a small ball $B_\epsilon(0)$ around $0$ and such that $X_\alpha
    (0) = e_\alpha$ for all $\alpha = 1, \ldots, 0$. It is well-known
    that such vector fields exist. Moreover, let $\Theta \in
    \Gamma^\infty(\Anti^2E|_U)$ be an arbitrary section of $\Anti^2
    E|_U$, locally written as $\Theta = \frac{1}{2}
    \Theta^{\alpha\beta} e_\alpha \wedge e_\beta$. Then $\theta =
    \frac{1}{2} \pi^*\Theta^{\alpha\beta} X_\alpha \wedge X_\beta \in
    \Xver^2(E|_U)$ is a Poisson structure such that $\theta_p$ has
    compact support around $0_p \in E_p$ \emph{and} $\theta (0_p) =
    \Theta^\ver(0_p)$ for all $p \in U$. Hence there are `many'
    vertical Poisson structures with compactly supported $\theta_p$.
    In the global situation we still have many of them but it is not
    clear whether we also can arrange to get every vertical lift at
    the zero section.
\end{example}

Let us now turn to star products. A star product $\star =
\sum_{r=0}^\infty \lambda^r C_r$ is called \emph{vertical} if all
$C_r$ are vertical bidifferential operators. In particular, the
antisymmetric part of $C_1$ defines a vertical Poisson structure
$\theta$ by $\{f, g\}_\theta = \frac{1}{\I} \left(C_1(f, g) - C_1(g,
    f)\right)$. In this case we say that $\star$ \emph{quantizes}
$\theta$.

A star product $\star$ is vertical iff
\begin{equation}
    \label{eq:StarVertical}
    f \star \pi^*u = f \pi^*u = \pi^*u \star f
\end{equation}
for all $f \in C^\infty(E)$ and $u \in C^\infty(M)$. Two vertical star
products $\star$ and $\star'$ are \emph{vertically equivalent} is
there exists a formal series $S = \id + \sum_{r=1}^\infty 
\lambda^r S_r$ of vertical differential operators $S_r$ such that
\begin{equation}
    \label{eq:StarEquivStar}
    S(f \star' g) = Sf \star Sg
    \quad
    \textrm{and}
    \quad
    S1 = 1.
\end{equation}
Analogously to vertical Poisson structures we can also restrict
vertical star products $\star$ to star products $\star_p$ for
$C^\infty(E_p)[[\lambda]]$. Clearly, $\star_p$ is still associative by
\eqref{eq:RestrictComposition} and we have
\begin{equation}
    \label{eq:RestrictionAlgMorph}
    \iota^*_p (f \star g) = \iota^*_p f \star_p \iota^*_p g.
\end{equation}
Moreover, vertically equivalent $\star$ and $\star'$ restrict to
equivalent $\star_p$ and $\star'_p$.

Using the vertical formality theorem one immediately obtains the
following existence and classification theorem by general arguments
analogous to \cite{kontsevich:2003a}:
\begin{theorem}[Vertical star products]
    \label{theorem:ExistenceClassification}
    Let $\pi: E \longrightarrow M$ be a vector bundle and $\uUver$ the
    vertical formality from Theorem~\ref{theorem:VerticalFormality}.
    \begin{enumerate}
    \item For a formal vertical Poisson structure $\theta \in
        \Xver^2(E)[[\lambda]]$ the definition
        \begin{equation}
            \label{eq:StarTheta}
            f \star_\theta g = fg + \sum_{r=1}^\infty
            \left(\frac{\I\lambda}{2}\right)^r
            \left(\Uver^{(r)}(\theta, \ldots, \theta)\right)(f, g)
        \end{equation}
        yields a vertical star product $\star_\theta$ quantizing
        $\theta_0$. If $\theta = \cc{\theta}$ is real,
        $\star_\theta$ is Hermitean.
    \item The map $\theta \mapsto \star_\theta$ induces a bijection on
        the level of vertical equivalence classes of formal vertical
        Poisson structures and vertical star products, respectively.
    \item The restriction $(\star_\theta)_p$ coincides with
        $\star_{\theta_p}$ which is the Kontsevich star product
        corresponding to $\theta_p$ on the vector space $E_p$ for $p
        \in M$.
    \end{enumerate}
\end{theorem}

%
% references
%

\begin{footnotesize}
%    \renewcommand{\arraystretch}{0.5} 
%    \bibliographystyle{ewde}
%    \bibliography{dqarticle,dqbook,dqprocentry,dqproceeding,preprints,script,dqthesis}

\end{footnotesize}

\end{document}